\documentclass[12pt, leqno]{article}
\usepackage{amsmath}
\usepackage{amssymb}
\usepackage{amsthm}
\usepackage[mathscr]{eucal}
\usepackage{latexsym}
\usepackage{amscd}
\usepackage{epic}
\usepackage{epsfig}
\usepackage{psfrag}

\begin{document}
\baselineskip=15pt

\voffset -1.5truecm
\oddsidemargin .5truecm
\evensidemargin .5truecm

\theoremstyle{plain}
\newtheorem{prop}{Proposition}
\newtheorem{coro}{Corolary}
\newtheorem{teor}{Theorem}
\newtheorem{demo}{Proof of Theorem 1}
\newtheorem{ejem}{Example}

\renewcommand{\refname}{\large\bf References}
\renewcommand{\thefootnote}{\fnsymbol{footnote}}

\def\flecha{\longrightarrow}
\def\asocia{\longmapsto}
\def\sii{\Longleftrightarrow}
\def\vector#1#2{({#1}_1,\dots,{#1}_{#2})}
\def\conjunto#1{\boldsymbol{[}\,{#1}\,\boldsymbol{]}}
\def\particion{\vdash}

\def\natural{{\mathbb N}}
\def\entero{{\mathbb Z}}
\def\racional{{\mathbb Q}}
\def\real{{\mathbb R}}
\def\complejo{{\mathbb C}}

\def\liri{Littlewood-Richardson\ }
\def\beze{Berenstein-Zelevinsky\ }
\def\cont{content\ }
\def\lr#1{{\sf LR}_{#1}}
\def\lrtipo#1{{\sf LR}_{#1}(\lambda,\mu,\nu)}
\def\pesos#1{{\sf D}_{#1}}
\def\triangulos#1{\real^{\binom{k+2}{2}-1}}
\def\tipo{(\lambda,\mu,\nu)}
\def\size#1{\vert {#1} \vert}
\def\arreglo#1{({#1}_{ij})_{\,0\le i\le j\le k}}
\def\hive#1{{\sf H}_{#1}}
\def\hivetipo#1{{\sf H}_{#1}(\lambda,\mu,\nu)}
\def\bz#1{{\sf BZ}_{#1}}
\def\bztipo#1{{\sf BZ}_{#1}(\lambda,\mu,\nu)}
\def\lrcoef{c^\lambda_{\mu\,\nu}}
\def\ent#1#2#3{{#1}_{{#2}\,{#3}}}

\def\carre#1{{\sf H}_{#1}(\lambda,\nu,\mu)}

\def\forma#1{{\sf sh}({#1})}
\def\sst{semistandard\ tableau\ }
\def\sucesion#1#2{{#1}_1 < {#1}_2 < \cdots < {#1}_{#2}}
\def\parti#1{\pi(#1)}

\def\HG{\Delta}
\def\BZ{\Gamma}
\def\vp{\varphi}


\setlength\unitlength{0.08em}
\savebox0{\rule[-2\unitlength]{0pt}{10\unitlength}%
\begin{picture}(10,10)(0,2)
\put(0,0){\line(0,1){10}}
\put(0,10){\line(1,0){10}}
\put(10,0){\line(0,1){10}}
\put(0,0){\line(1,0){10}}
\end{picture}}

\newlength\cellsize \setlength\cellsize{18\unitlength}
\savebox2{%
\begin{picture}(18,18)
\put(0,0){\line(1,0){18}}
\put(0,0){\line(0,1){18}}
\put(18,0){\line(0,1){18}}
\put(0,18){\line(1,0){18}}
\end{picture}}
\newcommand\cellify[1]{\def\thearg{#1}\def\nothing{}%
\ifx\thearg\nothing
\vrule width0pt height\cellsize depth0pt\else
\hbox to 0pt{\usebox2\hss}\fi%
\vbox to 18\unitlength{
\vss
\hbox to 18\unitlength{\hss$#1$\hss}
\vss}}
\newcommand\tableau[1]{\vtop{\let\\=\cr
\setlength\baselineskip{-16000pt}
\setlength\lineskiplimit{16000pt}
\setlength\lineskip{0pt}
\halign{&\cellify{##}\cr#1\crcr}}}
\savebox3{%
\begin{picture}(15,15)
\put(0,0){\line(1,0){15}}
\put(0,0){\line(0,1){15}}
\put(15,0){\line(0,1){15}}
\put(0,15){\line(1,0){15}}
\end{picture}}
\newcommand\expath[1]{%
\hbox to 0pt{\usebox3\hss}%
\vbox to 15\unitlength{
\vss
\hbox to 15\unitlength{\hss$#1$\hss}
\vss}}


\newcommand{\la}{\lambda}
\newcommand{\cd}{\mathcal D}
\newcommand{\cx}{\mathcal X}
\newcommand{\cy}{\mathcal Y}
\newcommand{\cp}{\mathcal P}
\newcommand{\cjp}{\mathcal{J}}
\newcommand{\nn}{\mathbb{N}}
\newcommand{\cR}{\mathcal R}
\newcommand{\cS}{\mathcal S}
\newcommand{\cQ}{\mathcal Q}

\newcommand{\CC}{\gimel}
\newcommand{\bs}{\mathbf{s}}

\newcommand{\wh}{\widehat}
\newcommand{\bu}{\bullet}
\newcommand{\emp}{\varnothing}
\newcommand{\sq}{\square}
\newcommand{\st}{\star}
\newcommand{\bd}{\blacktriangledown}
\newcommand{\ap}{\approx}
\newcommand{\aA}{\mathbb{A}}
\newcommand{\zz}{\mathbb{Z}}

\newcommand{\mto}{\mapsto}

\newcommand{\ff}{\Bbb F}
\newcommand{\cc}{\Bbb C}
\newcommand{\rr}{\Bbb R}
\newcommand{\qqq}{\Bbb Q}
\newcommand{\tT}{\Bbb T}
\newcommand{\ee}{\Bbb E}
\newcommand{\sS}{\Bbb S}
\newcommand{\oo}{\Bbb O}
\newcommand{\ooo}{\overline \oo}
\newcommand{\gG}{\Bbb G}
\newcommand{\sm}{\setminus}
\newcommand{\pr}{\prime}
\newcommand{\Ga}{\Gamma}
\newcommand{\Om}{\Omega}
\newcommand{\De}{\Delta}
\newcommand{\Si}{\Sigma}

\newcommand{\ga}{\gamma}
\newcommand{\si}{\sigma}
\newcommand{\de}{\delta}
\newcommand{\ep}{\epsilon}
\newcommand{\al}{\alpha}
\newcommand{\be}{\beta}
\newcommand{\om}{\omega}
\newcommand{\ka}{\kappa}
\newcommand{\ve}{\varepsilon}
\newcommand{\vr}{\text{\rm deg}}
\newcommand{\vs}{\varsigma}
\newcommand{\vk}{\varkappa}
\newcommand{\CR}{\mathcal R}
\newcommand{\GG}{\goth G}
\newcommand{\GI}{\goth I}
\newcommand{\CS}{\mathcal S}
\newcommand{\cC}{\mathcal C}
\newcommand{\ca}{\mathcal A}
\newcommand{\cb}{\mathcal B}
\newcommand{\CB}{\mathcal B}
\newcommand{\T}{\bold T}
\newcommand{\CT}{\mathcal T}
\newcommand{\BB}{\bold B}
\newcommand{\bP}{\bold P}
\newcommand{\ssu}{\subset}
\newcommand{\sss}{\supset}
\newcommand{\wt}{\widetilde}

\newcommand{\La}{\Lambda}
\newcommand{\const}{\text{{\rm const}}}
\newcommand{\rc}{ {\text {\rm rc}  } }
\newcommand{\y}{ {\text {\rm y}  } }
\newcommand{\sy}{ {\text {\rm sy}  } }
\newcommand{\Z}{ {\text {\rm Z} } }
\newcommand{\FF}{ {\text {\rm F} } }
\newcommand{\RR}{ {\text {\rm R} } }
\newcommand{\rk}{\dim}
\newcommand{\srk}{\dim}
\newcommand{\prev}{\prec}
\newcommand{\oa}{\overrightarrow}
\newcommand{\Ups}{\Upsilon}
\newcommand{\vt}{\vartheta}
\newcommand{\dv}{\divideontimes}
\newcommand{\cw}{\mathcal W}
\newcommand{\U}{{\text {\rm U} } }
\newcommand{\rT}{{\text {\rm T} } }
\newcommand{\Q}{{\text {\rm Q} } }
\newcommand{\mix}{{\text {\rm mix} } }
\newcommand{\sign}{{\text {\rm Sign} } }
\newcommand{\PSL}{{\text {\rm PSL} } }
\newcommand{\SL}{{\text {\rm SL} } }


\newcommand{\YT}{{\text {\rm YT} } }
\newcommand{\LR}{{\text {\rm LR} } }
\newcommand{\LRN}{{{\text {\textmd{LR}} } }}
\newcommand{\FV}{{\text {\rm CF} } }
\newcommand{\RSK}{{\text {\rm RSK} } }
\newcommand{\Mat}{{\text {\rm Mat} } }
\newcommand{\Can}{{\text {\rm Can} } }
\newcommand{\N}{{{\text {\textmd{N}} } }}

\newcommand{\ba}{\mathbf{a}}
\newcommand{\bb}{\mathbf{b}}
\newcommand{\bdd}{\mathbf{d}}
\newcommand{\bi}{\mathbf{i}}
\newcommand{\bm}{\mathbf{m}}

\newcommand{\bw}{\triangledown}

\newcommand{\red}{\hookrightarrow}
\newcommand{\att}{\star}
\newcommand{\di}{\diamond}
\newcommand{\vdi}{\diamond}
\newcommand{\bsl}{\langle}
\newcommand{\bsr}{\rangle}

\newcommand{\comp}{\circ}
\newcommand{\too}{\longrightarrow}
\newcommand{\lrm}{\hookrightarrow}

\newcommand{\weight}{{\text {\rm {\bf weight}} } }
\newcommand{\word}{{\text {\rm {\bf word}}}}

\newcommand{\ER}{\vskip1.cm
\hskip1.cm {\bf ****** \qquad !`
Ernesto, please write this part ! \qquad ******}
\vskip1.cm
}

\newcommand{\ERA}{\vskip1.cm
\hskip1.cm {\bf ****** \qquad !`
Ernesto, please add anything you like ! \qquad ******}
\vskip1.cm
}

\begin{centering}
{\Large\bf  Reductions of Young tableau bijections }\\[1.cm]
{{\large\sf  }}\vskip.1cm
{{\large{\bf Igor Pak}
\vskip.1cm
Massachusetts Institute of Technology\\
Cambridge, MA 02138 USA\\
e-mail: {\tt pak@math.mit.edu}\\}
\vskip.4cm
{\large\sf and}\vskip.4cm
{{\large{\bf Ernesto Vallejo}
\vskip.1cm
Instituto de Matem\'aticas, Unidad Morelia, UNAM\\
Apartado Postal 61-3, Xangari\\
58089 Morelia, Mich., MEXICO\\
e-mail: {\tt vallejo@matmor.unam.mx}}}
\vskip.6cm
June 23, 2004}
\\
\end{centering}

\vskip 2pc
\begin{abstract}
We introduce notions of linear reduction and linear equivalence
of bijections for the purposes of study bijections between Young
tableaux. Originating in Theoretical Computer Science, these notions
allow us to give a unified view of a number of classical
bijections, and establish formal connections between them.
\end{abstract}

\vskip 2.5pc

\section*{Introduction} \label{sec:intro}

Combinatorics of Young tableaux is a beautiful subject which
originated in the works of Alfred Young over a century ago,
and has been under intense development in the past decades
because of its numerous applications \cite{Fu,Ma,Sa,St}.
The amazing growth of the literature and a variety of
advanced extensions and generalizations led to some
confusion even about the classical combinatorial results in
the subject.  It seems that until now, researchers in the field
have not been able to unify the notation, techniques, and
systematize their accomplishments.

In this paper we concentrate on bijections between Young
tableaux, a classical and the most a combinatorially effective
part of subject.
The notable bijections include Robinson-Schensted-Knuth
correspondence, Jeu de Taquin, Sch\"utzenberger involution,
Littlewood-Robinson map, and
Benkart-Sottile-Stroomer's tableau switching.
The available descriptions of these bijections are so different,
that  a casual reader receives the impression that
all these maps are only vaguely related to each other.
Even though there is a number of well-known and
important connections between some of these bijections,
these results are sporadic and until this work did not fit
any general theory.
The idea of this paper is to give a formal general approach
to positively resolve this problem, and place these bijections
under one roof, so to speak.

We introduce new notions
of \emph{linear reduction} and \emph{linear equivalence} of
bijections, and show that the above mentioned bijections are
linearly equivalent.  This gives the first rigorous result
showing that these bijections are ``all the same'' in a certain
precise sense. A benefit of this approach is that we establish
a number of unexpected connections between classical
Young tableau bijections, and even discover a few
``traditional style'' conjectures.
We elaborate here on the nature and the origin of
our approach, while leaving definitions and main results
to the next two sections.

The philosophy behind this paper is the basic idea that Young tableau
bijections are significantly different from almost all classical bijections
in combinatorics, and that a complexity approach captures this difference.
On the one hand, the universe of combinatorial bijections is quite
heterogenous: some bijections are canonical while some are inherently
non-symmetric, some are recursive while some are direct and
explicit, some are very natural and easy to find while some
are highly nontrivial and their discovery is a testament to
human ingenuity (see e.g.~\cite{P1,St}).  On the other hand,
there is a property that almost all of them share and this
is that they can be computed in the linear number of steps (see below).
This is what explains why proofs of their validity tend to be relatively
straightforward, even if bijections' construction may seem complicated.
Similarly, the way bijections translate combinatorial statistics
between sets, tend to be relatively transparent,
which accounts for effectiveness of such bijections in the
majority of successful applications.

In contrast, Young tableau bijections are substantially
harder to establish, their working is delicate and proving their
properties is subtle (as is witnessed by the validity of the
Jeu de Taquin, the duality of the RSK correspondence, and the
many decades that passed before the
Littlewood-Robinson correspondence was formally proved).  To
explain this phenomenon we make two observations.  First, we
observe that all these bijections are inherently nonlinear and
require\footnote{We do not prove the lower bound, only the
upper bound on the time complexity.}
 roughly~$O(N^{3/2})$ number of arithmetic and logical
operations (in the size~$N$ it takes to encode the tableaux).
Second, we show that any of these bijections can be used to
construct any other.  This explains the identical power~$3/2$
in all cases, and at the same time asserts that all these
``exceptionally hard'' bijections are ``essentially the same''
and thus form a single class of ``counterexamples to the rule''
(that all bijections can be computed in linear time).
Making these claims rigorous is a difficult task which
required the introduction of new notation, definitions, and tools.

Following ideas of Theoretical Computer Science, we view
each \emph{bijection} as an algorithm with one type of combinatorial
objects as the input, and another as an output.  To make a
distinction, we say that a \emph{correspondence} is a one-to-one
map established (produced)
by a bijection.  Thus, several different bijections
can \emph{define} the same correspondence (cf.~\cite{P1}).
The \emph{complexity}, or the \emph{cost} of the algorithm, is,
roughly, the number of steps in the bijection.  One can
think of a correspondence  as a \emph{function} which
is computed by the algorithm (a bijection).  In contrast with
the emphasis in Cryptography, our correspondences
(and their inverses!) are relatively easy to compute.  On the
contrary, the bijections we analyze play an important role in
Algebraic Combinatorics in part due to this fact.

Now, as it is the case in Computational Complexity, finding
lower bounds for the cost of bijections defining
a given correspondence is a hard problem; we do not approach
this question in the paper.  Instead, we utilize what is known
as \emph{Relative Complexity}, an approach based on reduction
of combinatorial problems.  The reader may recall that the
class of \textsc{NP}-complete problems is closed under
polynomial time reductions~\cite{GJ}.  Similar notions exist
for a variety of problems in low complexity classes, with
various restrictions on time and space of the algorithms
(see e.g.~\cite{Pap}).

In this paper we consider only \textbf{linear time reductions};
since the bijections we consider require subquadratic time the
reductions have to preserve that.
Formally, we say that function~$f$ {\em reduces linearly} to~$g$,
if one can compute~$f$ in time linear in the time it takes
to compute~$g$.  We say
that~$f$ and~$g$ are \emph{linearly equivalent} if~$f$
reduces linearly to~$g$ and vice versa.  This defines
an equivalence relation on functions; it now can be translated
into a linear equivalence on combinatorial bijections.

Our main result is a linear equivalence of the Young tableau
bijections mentioned above, as well as few other known and new
bijections.  To present a (skew, semistandard) Young tableau
with~$\le k$ (possibly empty) rows and entries~$\le k$,
we need to write $\binom{k+1}{2}$
integers~$a_{i,j}$ which
represent the number of~$j$'s in~$i$-th row.  Ignoring the
size of~$a_{i,j}$, this gives input of size~$O(k^2)$.
Now, we shall see that all the bijections described above use the
same subquadratic number~$O(k^3)$ of arithmetic and logical
operations\footnote{Taking~$N =k^2$ this gives $O(N^{3/2})$
time mentioned above.}.
This is in sharp contrast with the Young tableau bijections
of linear cost, defined in the previous paper by the authors~\cite{PV}.
Roughly, we showed there that~$O(k^2)$ is the cost of bijections
between several combinatorial interpretations of the
Littlewood-Richardson's  coefficients.  This comes from
the fact that bijections in~\cite{PV} are given by simple
linear transformations, while Young tableau maps in this
paper are inherently nonlinear.


Despite our exhaustive search through the literature, it seems that
Computational Complexity has never been used in this area of
Algebraic Combinatorics.  In fact, we were able to find very
few references with any kind of computational approach
(see~\cite{BF,stem} for rare examples).  Of course, this is in
sharp contrast with other parts of Combinatorics such as
Graph Theory, Discrete Geometry, or Probabilistic Combinatorics,
where computational ideas led to important advances and
shift in philosophy.
%
%
We hope this paper will serve a
starting point in the future investigations of complexity
of combinatorial bijections.

\smallskip

We now elaborate on the content and the structure of the paper.
As the reader will see, this paper is far from being self-contained.
In fact, we never even define some of the classical Young tableau
bijections and in most proofs we assume that
the reader is familiar with the subject.
This decision was largely forced upon us, to keep the paper
of manageable size.
On the other hand, we are careful to include a number
of propositions giving properties and often alternative definitions
of these bijections.  Thus, much of the paper (the nontechnical part)
can be read by the reader unfamiliar with the subject, although
in this case some of our results may seem unmotivated, inelegant,
 or even simplistic.

The structure of the paper is unusual as we try
to emphasize the results themselves rather than the technical
details in the proofs.  We start by giving basic definitions and
listing the classical maps, in terms of combinatorial
objects they act upon (section~\ref{sec:basic}).
Even the reader well familiar with the subject is encouraged to
quickly go through these to become familiar with our notation.
We then present our new framework (of linear reductions)
and immediately state main
results (Section~\ref{sec:main}).  In Section~\ref{sec:prop}
we present properties of the bijections, connecting them to each
other through known results in the literature.
Section~\ref{sec:zoo} consists of a number of small subsections
which give linear reductions between various pairs of these
bijections.   This section represents the main part of the proof;
proofs of technical lemmas and other details are given in
Section~\ref{sec:proofs}.  In Section~\ref{sec:further}
we describe several conjectures and other important bijections
worthy of analysis.  Final remarks are given in
Section~\ref{sec:final}.

We should mention that throughout the first four sections we
make no references to the literature.  Instead, in
Subsection~\ref{sec:proofs-prop} we present a very brief
overview of the literature together with citation of sources
containing the propositions.  Further references are given
in Section~\ref{sec:further}.

\smallskip

{\bf Notation.}  \  We denote partitions by Greek letters:
$\la, \mu, \nu, \pi, \si, \tau,\ldots$ while maps are denoted
by different Greek letters:
$\vp, \psi, \phi, \zeta, \xi, \rho, \eta, \ldots$
Young tableaux are denoted by $A, B, C,\ldots$
generic sets of partitions are denoted by $\ca,\cb,\cC,\ldots$
and integer arrays (weights of the tableaux) are denoted
by $\ba,\bb,\bm,\ldots$  All Young diagrams and Young tableaux
are presented in English notation~\cite{Fu,St}.
Finally, we use $\nn = \{1,2,\ldots\}$ and
$\zz_{\ge 0} = \{0,1,2,\ldots\}$.

We should alert the reader to the fact that we use
\textit{``Proposition''}
mainly to describe known results, and we reserve \textit{``Theorem''}
and \textit{``Lemma''} for the new results.


\section{Basic definitions} \label{sec:basic}

\subsection{Young diagrams and Young tableaux}
\label{sec:basic-yd}

\begin{figure}[hbt]
\vskip.3cm
\psfrag{l}{$\la$}
\psfrag{m}{$\mu$}
\psfrag{n}{$\nu$}
\psfrag{p}{$\pi$}
\begin{center}
\epsfig{file=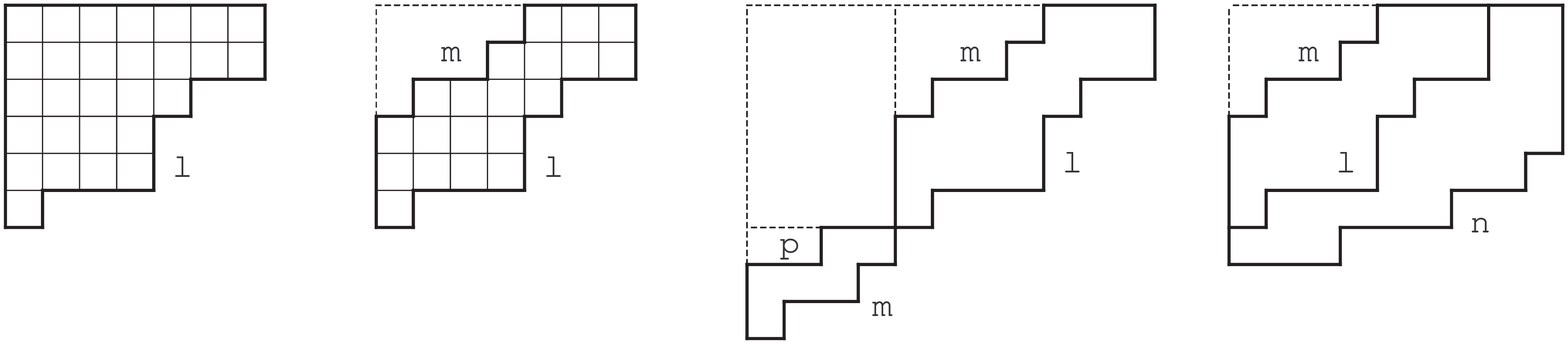,width=15cm}
\end{center}
\caption{Young diagram $[\la]$, and skew Young diagrams $[\la/\mu]$,
$[\mu/\pi]\comp [\la/\mu]$, and~$[\la/\mu] \star [\nu/\la]$, where
$\la = (7^254^21)$, $\mu=(431)$, $\nu = (9^4863)$, and~$\pi = (2)$. }
\label{f:defpar}
\end{figure}

A \emph{partition} $\la$ is a non-negative integer sequence
$(\la_1,\ldots,\la_\ell)$, such that
$\la_1\ge\ldots\ge\la_\ell\ge 0$\footnote{
Allowing $\la_i = 0$ is more natural from complexity point of
view and makes no difference for Young tableau bijections.}.
Denote by~$\ell = \ell(\la)$ the number of parts of~$\la$,
and let $|\la| = \la_1 + \ldots + \la_\ell$\,.
We say that~$\la$
is a \emph{partition of $n=|\la|$} and write $\la \vdash n$.
We represent a partition
graphically by a \emph{Young diagram}~$[\la]$ defined to be
a collection of squares
$\{(i,j) \in \zz^2 \mid 1 \le j \le \la_i, \, 1 \le i \le \ell \}$
(see Figure~\ref{f:defpar}).
Throughout the paper we often make no distinction
between partitions and the corresponding Young diagrams.

We say that $\mu \ssu \la$ if $\mu_i \le \la_i$ for all~$i>0$.
In other words, we have $\mu \ssu \la$ for partitions whenever
$[\mu] \ssu [\la]$ for Young diagrams viewed as sets of squares
(see Figure~\ref{f:defpar}).
%
%
A \emph{skew Young diagram} $[\la/\mu]$ is the shape of a
set of squares in~$[\la] \,- \,[\mu]$, where~$\mu \ssu \la$.
Let~$|\la/\mu| = |\la| - |\mu|$ denote the number of squares
in~$[\la/\mu]$, and $\ell(\lambda/\mu)$ the height of $[\la/\mu]$.
Without loss of generality we can always assume that $\ell(\lambda/\mu)=
\ell(\lambda)$.
We say that a skew Young diagram~$[\la/\mu]$ is \emph{attached} to
a skew Young diagram~$[\mu/\nu]$, for all~$\nu \ssu \mu \ssu \la$,
and denote by~$[\la/\nu] = [\mu/\nu] \att [\la/\mu]$ the union
of these two diagrams.
For two skew diagrams $[\la/\mu]$ and
$[\nu/\tau]$ we define a \emph{composition}
\[
[\nu/\tau] \comp  [\la/\mu] :=
\bigl[\,(\la_1+\nu_1,\dots,\la_\ell+\nu_1,\nu_1,\dots,\nu_r)/
(\mu_1+\nu_1,\dots,\mu_\ell+\nu_1,\tau_1,\dots,\tau_r)\,\bigr],
\]
where $\ell = \ell(\la)$, $r = \ell(\nu)$.
Graphically, this corresponds to placing $[\la/\mu]$ above and
to the right of $[\nu/\tau]$ (see Figure~\ref{f:defpar}).
We can generalize this to $[\nu/\tau] \comp_{a,b}  [\la/\mu]$, where
$[\la/\mu]$ is shifted by~$b$ squares above $[\nu]$ and by~$a$
squares to the right of the left margin of~$[\nu]$, here $a\ge \nu_1$.
Thus, $[\nu/\tau] \comp_{\nu_1,0}  [\la/\mu] = [\nu/\tau] \comp  [\la/\mu]$.


\begin{figure}[hbt]
\vskip.3cm
\psfrag{A}{$A$}
\psfrag{B}{$B$}
\psfrag{C}{$C$}
\psfrag{l}{$\la$}
\psfrag{m}{$\mu$}
\psfrag{n}{$\nu$}
\begin{center}
\epsfig{file=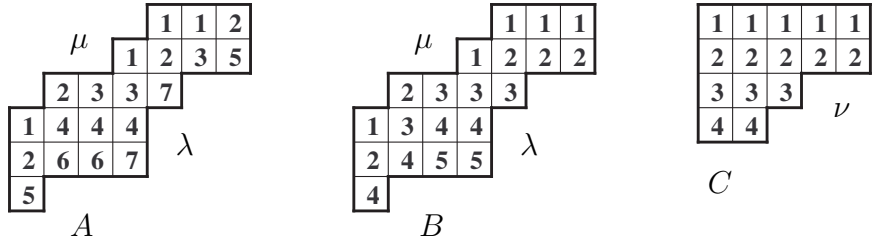,width=12.cm}
\end{center}
\caption{Young tableaux $A \in \YT(\la/\mu,\tau)$,
$B \in \LR(\la/\mu,\si)$, and $C = \Can(\nu)$, where
$\la = (7,7,5,4,4,1)$, $\mu=(4,3,1)$,
$\tau = (4,4,3,3,2,2,2)$, $\si=(5,5,4,4,2)$, and $\nu = (5,5,3,2)$. }
\label{f:deftab}
\end{figure}

A \emph{Young tableau}~$A$ of shape~$\la/\mu$ is a
function $A : [\la/\mu] \to \nn$, which is increasing in columns
and non-decreasing in rows.   We think of the values $A(i,j)$ as
being written in the squares of the (skew) Young diagram
(see Figure~\ref{f:deftab}).  The \emph{weight} of a tableau~$A$,
denoted $\weight (A)$, is a sequence $\bm = (m_1,m_2,\dots)$, where
$m_r = \bigl|\{(i,j) \in [\la/\mu] \mid A(i,j) = r\}\bigr|$.
Clearly, $m_1+m_2+\ldots = n = |\la /\mu|$.
Denote by~$\ell(\bm)$ the length of~$\bm$ and by
$\YT(\la,\bm)$ the set of Young tableaux of shape~$\la$
and weight~$\bm$.
In some cases the weight~$\bm$ will be a partition of~$n$.
Let $\YT(\la/\mu;  k)$ denote the set of Young tableaux of
shape $\la/\mu$ and weight~$\bm$, such that~$\ell(\bm)\le k$.
In other words, $\YT(\la/\mu; \, k)$ contains all tableaux
with integers~$\le k$.
By $\Can(\la)$ we
denote the unique Young tableau~$A \in \YT(\la,\la)$ with
1's in the first row, 2's in the second row, etc.
We call such~$A$ a \emph{canonical tableau} of shape~$\la$
(see Figure~\ref{f:deftab}).

Denote by $A\comp B$ the natural composition of two (skew)
Young tableaux. Similarly,
denote by $A\att B$ the attachment of two (skew)
Young tableaux, whenever this is possible.

For a sequence of integers~$\bi = (i_1,\ldots,i_n)$, $i_r \in \nn$,
denote by $m_j(\bi)$ the number of $j$'s in~$\bi$.
We say that~$\bi$ is \emph{positive} if
$m_1(\bi) \ge m_2(\bi) \ge \cdots$.
Similarly, we say
that~$\bi = (i_1,\ldots,i_n)$ is \emph{dominant}
if~$\bi_r = (i_1,\ldots,i_r)$ is positive for all~$1 \le r \le n$.
For a Young tableau~$A$ denote by~$\word(A)$ the sequence obtained
by reading right-to-left the first row, then the second row, etc.

We say that~$A$ is a \emph{Littlewood-Richardson} (LR) \emph{tableau}
if~$A$ is a Young tableau and~$\word(A)$ is dominant.
Denote by~$\LR(\la/\mu,\nu)$ the set of all $\LR$-tableaux
of shape $\la/\mu$ and weight~$\nu$.  Observe that
when~$\mu = \emp$, there exist only one $\LR$-tableau:
a canonical tableau of shape~$\la$ and weight~$\nu=\la$.

\begin{figure}[hbt]
\vskip.3cm
\psfrag{A}{$A$}
\psfrag{B}{$B$}
\psfrag{C}{$C$}
\psfrag{l}{$\la$}
\psfrag{m}{$\mu$}
\psfrag{n}{$\nu$}
\psfrag{s}{$\si$}
\psfrag{t}{$\tau$}
\begin{center}
\epsfig{file=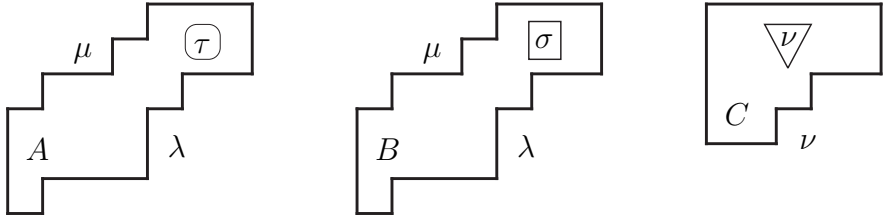,width=12.cm}
\end{center}
\caption{Illustrations of Young tableaux $A \in \YT(\la/\mu,\tau)$,
$B \in \LR(\la/\mu,\si)$, and $C = \Can(\nu)$. }
\label{f:defpict}
\end{figure}

To simplify the examples we illustrate the tableaux by
diagram drawings.  We circle the weight of general Young
tableaux, put it in a square in case of LR-tableaux,
and put it in a downward triangle for canonical tableaux.
See Figure~\ref{f:defpict} for an illustration of Young
tableaux given in Figure~\ref{f:deftab}.

\subsection{Main bijections}
\label{sec:basic-bij}

Here we present only the notation and set-theoretic statements
of the bijections.  These bijections will be further studied
in Section~\ref{sec:prop}.

\smallskip

Let $\ba = (a_1,\dots,a_k)$, $\bb = (b_1,\dots,b_k)$, such
that $a_i, b_j \in \zz_{\ge 0}$ and $|\ba| = |\bb|$, i.e.
$a_1+\cdots + a_k = b_1 + \cdots + b_k = N$.
Denote by $\Mat(\ba,\bb)$ the set of $k\times k$
matrices $V = (v_{i,j})$,
such that $v_{i,j} \in \zz_{\ge 0}$ \, and
$$\sum_{j=1}^k \, v_{i,j} \,  = \, a_i\,, \ \
\sum_{i=1}^k \, v_{i,j} \, = \, b_j\,, \ \ \text{for all} \ \, 1\le i,
\,j \le k\,.
$$

\textbf{1)} \,
The \emph{Robinson-Schensted-Knuth} (RSK) {\em correspondence} is a
one-to-one correspondence
$$\vp: \, \Mat(\ba,\bb) \, \too \, \bigcup_{\la \vdash N} \,
\YT(\la,\bb) \times \YT(\la,\ba),
$$
where $|\ba| = |\bb | = N$ as above.
Recall that if $\vp(V)=(B,A)$, then
$B$ is called the {\em insertion tableau} and $A$ is called the
{\em recording tableau}, that is, $B$ is the $P$-tableau and
$A$ is the $Q$-tableau.


\smallskip

\textbf{2)} \,
Let $\mu \ssu \la$ and $\ba=(a_1,\dots,a_m)$ be a sequence of
non-negative integers such that $n=|\ba|=|\la/\mu|$.
\emph{Jeu de Taquin} is a map
$$\psi: \, \YT(\la/\mu,\ba) \, \too \,
\bigcup_{\pi \vdash n} \, \YT(\pi,\ba).$$
We should emphasize that~$\psi$ is neither an into nor an
onto map.  As we shall see later, it plays a central role
in the universe of Young tableau bijections.

\smallskip

\textbf{3)} \,
In notation of the Jeu de Taquin map, the
\emph{Littlewood-Robinson map} is given by the following
one-to-one correspondence:
$$\phi: \, \YT(\la/\mu,\ba) \, \too \,
\bigcup_{\pi \vdash n} \, \YT(\pi,\ba) \times
\LR(\la/\mu,\pi).
$$


\smallskip

\textbf{4)} \,
Suppose now that $\mu \ssu \la$ and ${\mathbf c}$,
${\mathbf d}$ are sequences of non-negative integers
such that $n = |\la/\mu| = |{\mathbf c}| + |{\mathbf d}|$.
\emph{Tableau Switching} is a one-to-one correspondence:
\[
\zeta: \, \bigcup_{\pi \, \vdash |\la|-|{\mathbf c}|} \,
\YT(\pi/\mu,{\mathbf d}) \times \YT(\la/\pi,{\mathbf c})  \ \too
\bigcup_{\si \, \vdash |\la|-|{\mathbf d}|} \,
\YT(\si/\mu,{\mathbf c}) \times \YT(\la/\si,{\mathbf d})  .
\]


\textbf{4$^\mathbf{\prime}$)} \,
We shall need a special notation for Tableau Switching in
case $\mu$ is the empty partition, that is for {\em normal shapes}:
$$\zeta^\N: \, \bigcup_{\pi \, \vdash |{\mathbf d}|} \,
\YT(\pi,{\mathbf d}) \times \YT(\la/\pi,{\mathbf c})  \ \too
\bigcup_{\si \, \vdash |{\mathbf c}|} \,
\YT(\si,{\mathbf c}) \times \YT(\la/\si,{\mathbf d}) .$$


\smallskip

\textbf{5)} \,
Let $\ba = (a_1,\dots,a_m)$ be an integer sequence, such that
$|\ba| = |\la/\mu|$.
Consider a sequence~$\ba^\ast := (a_m,\dots,a_1)$.  Clearly,
$|\ba^\ast| = |\la/\mu|$ as well.  The \emph{Sch\" utzenberger
involution} is a one-to-one correspondence:
\[
\xi: \, \YT(\la/\mu,\ba) \, \to \, \YT(\la/\mu,\ba^\ast).
\]

\textbf{5$^\mathbf{\prime}$)} \,
We shall also need a special notation for the Sch\"utzenberger involution
for {\em normal shapes}:
$$\xi^\N: \, \YT(\la,\ba) \, \to \, \YT(\la,\ba^\ast).$$

\smallskip

\textbf{6)} \,
There is a different one-to-one correspondence of the same kind
as $\xi$ called {\em reversal}:
\[
\chi : \, \YT(\la/\mu,\ba) \, \to \, \YT(\la/\mu,\ba^\ast).
\]
This correspondence has certain advantages over the
Sch\"utzenberger involution, and will be useful for the proofs.

\smallskip

\textbf{7),~8)}  \,
Finally, suppose $\nu \vdash |\la/\mu|$. The \emph{Fundamental
Symmetry} is a one-to-one correspondence:
$$\rho: \, \LR(\la/\mu,\nu) \, \to \, \LR(\la/\nu,\mu).
$$
In this paper we define several versions of the
Fundamental Symmetry map, and the main results cover
the {\em First} and the {\em Second  Fundamental Symmetry}
maps.  In general, we say that a map \emph{gives the Fundamental
Symmetry map}
if it is a one-to-one correspondence between sets as above.

\bigskip

\subsection{Computational preliminaries}
\label{sec:basic-comp}

Let $D = (d_1, \dots, d_n)$ be an array of integers, and let
$m = m(D) := \max_i d_i$.
The \emph{bit-size} of~$D$, denoted by~$\bsl D \bsr$,
is the amount of space required to store~$D$.  For simplicity
everywhere below we assume that
$\bsl D \bsr = n \, \lceil\log_2 m+1 \rceil$.


We view a bijection $\tau: \ca \to \cb$ as an algorithm which
inputs $A\in \ca$ and outputs~$B=\tau(A) \in \cb$.  We need to
present Young tableaux as arrays of integers so that we can
store them and compute their bit-size.


Suppose $A \in \YT(\la/\mu; k)$.
An important way to encode $A$ is through a matrix, often
called the {\em Gelfand-Tsetlin} (GT) {\em pattern}
$(a_{i,j})$ of the tableau, whose entries satisfy certain
inequalities which are irrelevant for the purposes of this paper.
It is defined by $a_{i, 0} = \mu_i$ and $a_{i, j} =\mu_i~+$~the
number of integers in row $i$ which are $\le j$, for
$1\le i\le \ell(\la/\mu)$ and $1\le j\le k$.
The tableau~$A$ can be now be viewed as a matrix~$(a_{i, j})$;
this is the way Young tableaux will be presented in the input and
output of the algorithms.
Another useful way to encode $A$ is through its
{\em recording matrix} $(c_{i, j})$, which is defined by
$c_{i, j} = a_{i, j} - a_{i, j-1}$; in other words,
$c_{i, j}$ is the number of $j$'s in the $i$-th row of $A$.

Finally, we say that a map $\ga : \ca \to \cb$ is
\emph{size-neutral} if the ratio~$\bsl\ga(A)\bsr/\bsl A\bsr$ is
bounded for all $A \in \ca$.  Throughout the paper
we consider only size-neutral maps, often
without emphasizing it.  See the remark following Theorem~2
(Section~\ref{sec:main-results}) on the reasoning behind this
definition.

\vskip1.cm

\section{Reductions of bijections and main results}
\label{sec:main}

\subsection{Linear reductions}
\label{sec:main-red}

Think of a \emph{bijection}, or any \emph{explicit map}
in general, as an algorithm with input and output
written as an array of integers.  Hereafter by
\emph{size} of the input/output we mean bit-size.
As before, let~$\bsl A\bsr$ denote the bit-size of the
integer array~$A$.  We also say that a bijection or an
explicit map \emph{defines} a correspondence between input
and output set. Clearly, many different bijections can
define the same one-to-one correspondence.

Let $\ca$ and $\cb$ be two possibly infinite sets of
finite integer arrays, and let $\de : \ca \to \cb$
be an explicit map between them. We say that~$\de$ has
\emph{linear cost} if~$\de$ computes
$\de(A) \in \cb$ in \emph{linear time}~$O(\bsl A\bsr)$ for
all~$A \in \ca$.

\smallskip

There are many ways to construct new bijections out
of existing ones.  We call such algorithms \emph{circuits}
and define below several of them that we need.

\smallskip

$1)$ \ Suppose $\de_1: \ca \to \cx_1$,
$\ga: \cx_1 \to \cx_2$ and $\de_2: \cx_2 \to \cb$, such
that $\de_1$ and $\de_2$ have linear cost.
Consider $\chi = \de_2 \circ \ga \circ \de_1: \ca \to \cb$.
We call this circuit \emph{trivial} and denote it by
$I(\de_1,\ga,\de_2)$.

\smallskip

$2)$ \ Suppose $\ga_1: \ca \to \cx$,
$\ga_2: \cx \to \cb$, and let
$\chi = \ga_2 \circ \ga_1: \ca \to \cb$.
We call this circuit \emph{sequential} and denote it by
$S(\ga_1,\ga_2)$.

\smallskip

$3)$ \ Suppose $\de_1: \ca \to \cx_1 \times \cx_2$,
$\ga_1: \cx_1 \to \cy_1$, $\ga_2: \cx_2 \to \cy_2$,
and $\de_2: \cy_1 \times \cy_2 \to \cb$, such
that that $\de_1$ and $\de_2$ have linear cost.
Consider $\chi = \de_2 \circ (\ga_1 \times \ga_2)
\circ \de_1: \ca \to \cb$. We call
this circuit \emph{parallel} and denote it by
$P(\de_1,\ga_1,\ga_2,\de_2)$.

\medskip

\noindent
Fix a bijection~$\be$.
We say that~$\CC$ is a
\emph{$\be$-based ps-circuit} if one of the following holds:

\smallskip

$\bu$ \ $\CC = \de$, where $\de$ has linear cost.

$\bu$ \ $\CC = I(\de_1,\be,\de_2)$.

$\bu$ \ $\CC = P(\de_1,\ga_1,\ga_2,\de_2)$, where
$\ga_1,\ga_2$ are $\be$-based ps-circuits.

$\bu$ \ $\CC = S(\ga_1,\ga_2)$, where $\ga_1,\ga_2$ are
$\be$-based ps-circuits.

\smallskip

\noindent
In other words, $\CC$ is a $\be$-based ps-circuit if there
is a parallel-sequential algorithm which uses only a finite
number of linear cost maps and a finite number of maps~$\be$.
The $\be$-{\em cost} of $\CC$ is the number of times the map $\be$
is used; we denote it by $\bs(\CC)$.

Let $\ga:\ca \to \cb$ be a map produced by the
$\be$-based ps-circuit~$\CC$.  We say that $\CC$
\emph{computes~$\ga$ at cost~$\bs(\CC)$ of~$\be$}.

\medskip

We say that a map~$\al$ is \emph{linearly reducible}
to~$\be$, write $\al \lrm \be$,
if there exist a finite $\beta$-based ps-circuit~$\CC$ which
computes~$\al$.  In this case we say that~$\al$ \emph{can be
computed in at most~$\bs(\CC)$ cost of~$\be$.}

We say that maps~$\al$ and~$\be$ are
\emph{linearly equivalent}, write $\al \sim \be$,
if~$\al$ is linearly reducible
to~$\be$, and~$\be$ is linearly reducible to~$\al$.

\subsection{Main results}
\label{sec:main-results}

Our first result is the linear equivalence of Young tableau
bijections given in Section~\ref{sec:basic-bij}.


\bigskip

\noindent
{\bf Theorem~1.} \ {\it
The following maps are linearly equivalent:

\smallskip

$1)$ \ RSK map~$\vp$.

$2)$ \ Jeu de Taquin map~$\psi$.

$3)$ \ Littlewood-Robinson map~$\phi$.

$4)$ \ Tableau Switching map~$\zeta$.

$5)$ \ Sch\"utzenberger involution for normal shapes $\xi^\N$.

$6)$ \ Reversal $\chi$.

$7)$ \ First Fundamental Symmetry map~$\rho_1$.

$8)$ \ Second Fundamental Symmetry map~$\rho_2$.

\smallskip

\noindent
Moreover, each of these maps
can be computed in at most~$36$ times the cost of any other map.
}

\bigskip

The following theorem gives a positive result about the
efficient computation of maps~$1)-8)$.

\bigskip


\noindent
{\bf Theorem~2.} \ {\it For the eight maps $1)-8)$ as in
Theorem~1, let~$k$ and~$m$ be defined as follows:

%

\smallskip
\begin{tabular}{l l l}
$1)$ &  $k := \max\{\ell(\ba), \ell(\bb)\}$, \ &
$m := \max\{ \sum_i a_i, \sum_j b_j \}$, \\
$2)\,, 3)\,, 5)\,, 6)$ \  &  $k := \ell(\ba)$, &
$m := \la_1$, \\
$4)$ & $k := \max\{ \ell({\mathbf c}), \ell({\mathbf d}) \}$, &
$m := \la_1$, \\
$7)\,, 8)$ &  $k := \ell(\nu)$, &
$m := \la_1$.
\end{tabular}

\smallskip

\noindent
Then the image of maps $1)-8)$ can be computed at a
cost of~$O(k^3 \log m)$.
}

\bigskip

We should emphasize that the standard definitions of these maps
give a weaker result.  For example, Jeu de Taquin as defined
in the literature (see e.g.~\cite{Sa,St})
requires $O(|\mu|\cdot |\la/\mu|)$ square moves,
much greater than the bound above.

\medskip

\noindent
 {\bf Remark 1.} \
It may seem that by Theorem~1, it suffices to establish
the efficient computation of any one of the maps.  This
is not the case since we compare the maps by the number
of times other maps are used, not by the timing.
A priori it can (and does) happen that
the maps increase the bit-size of combinatorial
objects, when they transform the input into the output.
This affects the timing of the subsequent applications
of these maps.  To control this, we show that all maps
we consider are in fact size-neutral,
so Theorem~1 remains applicable in this case.
We formalize this observation in
Section~\ref{sec:proofs-thms}; for now we
suggest the reader simply views linear reductions as
if they were reductions on the time complexity of the
maps.


\section{Properties of bijections}
\label{sec:prop}

\subsection{Bender-Knuth transformations}
\label{sec:prop-bk}

Define \emph{Bender-Knuth} (BK) \emph{transformations}
$s_1, s_2, \dots$ as follows.
Consider a Young tableau~$A \in \YT(\la/\mu, \ba)$.
Let~$m=\ell(\ba)$ be
the length of~$\ba$, and let $(a_{i,\, j})$ be the corresponding
GT-pattern.  For any~$1\le r < m$, let~$B=s_r(A)$ be a Young
tableau such that the corresponding GT-pattern $(b_{i,\, j})$
is defined as follows:
$$
b_{i,\, j} \, = \, \left\{
\aligned
& \min\{a_{i,r+1},a_{i-1,r-1}\} + \max\{ a_{i,r-1},a_{i+1,r+1}\}
 - a_{i,r}
\,, \, \ \text{if} \ j=r, \\
& \, a_{i,j}\,, \, \ \text{otherwise.}
\endaligned
\right.
$$
Observe that BK-transformations are defined by bijections
$$
s_i : \, \YT(\la/\mu, \ba) \to \YT(\la/\mu, \ba^\pr), \ \ \ \
s_i : \, A \mto B,
$$
where $\ba = (a_1,\ldots,a_i,a_{i+1},\ldots,a_m)$, and
$\ba^\pr = (a_1,\ldots,a_{i+1},a_i,\ldots,a_m)$.  They also
satisfy the following relations:
\begin{equation}
s_i^2 = 1\,, \ \ \ s_i s_j = s_j s_i\,, \ \,
\text{if} \ \, |i-j| \ge 2\,.
\tag{$\lozenge$}
\end{equation}
Now define elements $t_{r,\, m-r}$ and $z_m$ as follows:
\[
z_m \, = \, (s_1) (s_2 s_1) (s_3 s_2 s_1) \cdots
(s_{m-1} \cdots s_2 s_1),
\]
\[
t_{r, \, m-r} \, = \, (s_{m-r} s_{m-r+1} \cdots s_{m-1}) \cdots
(s_{2} s_{3} \cdots s_{r+1} ) (s_{1} s_2 \cdots s_{r} ),
\]
where $1 \le r < m$.
By definition,
\[
z_m : \, \YT(\la/\mu, \ba) \to \YT(\la/\mu, \ba^\ast), \qquad
t_{r,\, m-r} : \, \YT(\la/\mu, \ba) \to \YT(\la/\mu, \bb),
\]
where $\ba^\ast = (a_m,\ldots,a_{r+1},a_{r},\ldots,a_1)$,
and
$\bb = (a_{r+1},\ldots,a_m, a_1\ldots,a_{r})$.
The maps satisfy the following well-known relations:
\begin{equation}
z_m^2\, = \, 1, \ \ \ \ \ t_{r,\, m-r} \, t_{m-r,\, r} \, = \, 1.
\tag{$\maltese$}
\end{equation}
We shall also need the following relation:
\begin{equation}
z_{l+k} \, = \, z_k \,  t_{l, \, k} \,  z_l\, .
\tag{$\circledast$}
\end{equation}
It will be proved in Section~\ref{sec:proof-rel}.
\bigskip

\begin{prop}\label{evac}
The \emph{Sch\" utzenberger involution map}~$\xi$
coincides with the map~$z_m$ defined as above: \, $\xi = z_m$.
\end{prop}

\medskip
Let $\mu \ssu \pi \ssu \la$ be partitions and ${\mathbf c}$, ${\mathbf d}$
be vectors of non-negative integers such that
$|\la / \mu|= |{\mathbf c}| + |{\mathbf d}|$.
Denote $r=\ell({\mathbf d})$ and $m= \ell({\mathbf c}) + \ell({\mathbf d})$.
Let $A\in \YT(\la/\pi,{\mathbf c})$, $B\in \YT(\pi/\mu, {\mathbf d})$.
Let $\widetilde{A}$ denote the tableau obtained from $A$ by adding $r$
to each entry of $A$.
Denote ${\mathbf e}=(c_1,\dots,c_{m-r},\, d_1,\dots, d_r)$ and
${\mathbf f}=(d_1,\dots,d_r,\, c_1,\dots, c_{m-r})$.
Consider $B\star \widetilde{A}\in \YT(\la/\mu, {\mathbf f})$;
thus, $t_{r,m-r}(B\star \widetilde{A})\in \YT(\la/\mu, {\mathbf e})$.
Decompose this tableau as $A^\prime \star \widetilde{B}$ with
$A^\prime \in \YT(\si/\mu, {\mathbf c})$, for some partition $\sigma$
of $|\mu| + |{\mathbf c}|$, and $\widetilde{B} \in
\YT(\la/\si)$.
Let $B^\prime$ be obtained from $\widetilde{B}$ by subtracting $m-r$
to each one of its entries, so that $B^\prime\in\YT(\la/\si, {\mathbf d})$.
Finally, let $\tilde{t}_{r,m-r}(B,A) = (A^\prime, B^\prime)$.


\begin{prop}\label{switch}
The \emph{Tableau Switching map}~$\zeta$ coincides with
the map~$\tilde{t}_{r,m-r}$ defined as above.
\end{prop}

\subsection{Jeu de Taquin is everywhere}
\label{sec:prop-jdt}

Consider the RSK map first. Let $\vp: V \mapsto (B,A)$, where
$V = (v_{i,j}) \in \Mat(\ba,\bb)$, and $\ell(\ba)=\ell(\bb)=k$.
Define $[\pi/\sigma] = [a_1] \circ [a_2] \circ \cdots \circ [a_k]$ and
$[\rho/\tau] = [b_1] \circ [b_2] \circ \cdots \circ [b_k]$.

\begin{prop}\label{rsk}
Let $Y$ be the Young tableau in $YT(\pi/\sigma, \bb)$ that has exactly
$v_{i,j}$ entries equal to $j$ in row $(k+1)-i$.
Then, for the \emph{Jeu de Taquin map} $\psi$, one has $\psi(Y)=B$.
Similarly, if $X$ is the Young tableau in $YT(\rho/\tau, \ba)$ that
has exactly $v_{j,i}$ entries equal to $j$ in row $(k+1)-i$.
Then $\psi(X)=A$.
\end{prop}


Let $A \in \YT(\la/\pi,\ba)$, $B \in \YT(\pi,\bb)$.
Denote by $(A^\pr,B^\pr) = \zeta^\N(B , A)$ their
tableau switching, where $A^\pr \in \YT(\si,\ba)$,
$B^\pr \in \YT(\la/\si,\bb)$ for some $\si \vdash |\la/\pi|$.

\begin{prop}\label{swijeu}
In the notation above,
the image $(A^\pr,B^\pr)$ of $(B,A)$ under the
\emph{Tableau Switching map}~$\zeta$ satisfies $A^\pr = \psi(A)$.
\end{prop}


One can obtain the Jeu de Taquin map~$\psi$ as a projection of the
Littlewood-Robinson map~$\phi$ onto the first component.

\begin{prop} \label{lirojeu}
Let the image of the \emph{Littlewood-Robinson map}~$\phi(A)= (B,C)$.
Then $\psi(A) = B$.
\end{prop}

Let $[\la]$ be  Young diagram, and let
$[r^\ell]$ be the smallest size rectangle
containing~$[\la]$, i.e.~$\ell = \ell(\la)$ and~$r = \la_1$.
Denote by $[\la^\bu]$ a skew Young diagram
$[r^\ell/(r-\la_\ell,\dots,r-\la_1)]$.  If
$A \in \YT(\la,\ba)$ is a Young tableau corresponding to
a recording matrix $C=(c_{i,\, j})$, let
$A^\bu \in \YT(\la^\bu,\ba^\ast)$ by a tableau corresponding to
a recording matrix $C^{\bu} = (c_{\ell+1-i,\, k+1-j})$.

\begin{prop}\label{evajeu}
The image of the \emph{Sch\"utzenberger involution map}~{\rm $\xi^\N(A)$}
coincides with the image~$\psi(A^\bu)$.
\end{prop}

Note that Propositions~\ref{swijeu} and~\ref{evajeu} do not
have natural extensions to skew Young diagrams.

\bigskip

\subsection{Hidden symmetries of RSK}
\label{sec:prop-RSK}

The first hidden symmetry is called \emph{duality} and
has already been used in Proposition~\ref{rsk}.

\begin{prop}
Suppose $\vp(V) = (B,A)$.
Then $\vp(V^\pr) = (A,B)$, where $V^\pr$ is the transpose matrix of $V$.
\label{rsksim}
\end{prop}

The second symmetry is as classical as duality but not as well-known.
For any $V=(v_{i,j})\in \Mat(\ba,\bb)$ of size~$k\times k$
we denote $V^* := (v_{k+1-i,\, k+1-j})\in \Mat(\ba^*,\bb^*)$.



\begin{prop}\label{rskeva}
In the notation above, suppose $\vp(V) = (B,A)$.
Then $\vp(V^*) = (\xi(B), \xi(A))$.
\end{prop}

\bigskip

\subsection{First Fundamental Symmetry map}
\label{sec:prop-lr}

We start with the following important characterization
of Littlewood-Richardson (LR) tableaux:

\begin{prop}\label{jeucan}
Suppose $A \in \YT(\la/\mu,\nu)$.
Then~$A$ is a LR-tableau if and only if $\psi(A)$ is the
canonical tableau~$\Can(\nu)$.
\end{prop}


Now we are ready to define
the \emph{First Fundamental Symmetry map}~$\rho_1$.
Let $A \in \YT(\la/\mu,\nu)$, and let $B = \Can(\mu)$.
Consider $(A^\pr,B^\pr) = \zeta(B,A)$, where
$A^\pr \in \YT(\si,\nu)$,
$B^\pr \in \YT(\la/\si,\mu)$ for some $\si \vdash |\la/\mu|$.
Define~$\rho_1(A) = B^\pr$.
By Propositions~\ref{swijeu} and~\ref{jeucan}, if
$A \in \LR(\la/\mu,\nu)$, we have~$A^\pr = \Can(\si)$, and
therefore~$\nu=\si$.
Similarly, since $B= \Can(\mu)$, by the involution property~$(\maltese)$
of the tableau switching, we have $B^\prime\in \LR(\la/\nu,\mu)$
and~$\rho_1^2 = 1$.

\begin{prop}
The map~$\rho_1$ is well defined
and gives a \emph{fundamental symmetry map}.
\label{funsim1}
\end{prop}

\begin{figure}[hbt]
\vskip.3cm
\psfrag{A}{$A$}
\psfrag{B}{$B$}
\psfrag{C}{$C$}
\psfrag{A'}{$A'$}
\psfrag{B'}{$B'$}
\psfrag{C'}{$C'$}
\psfrag{l}{$\la$}
\psfrag{m}{$\mu$}
\psfrag{n}{$\nu$}
\psfrag{s}{$\si$}
\psfrag{t}{$\tau$}
\psfrag{r}{$\rho_1$}
\begin{center}
\epsfig{file=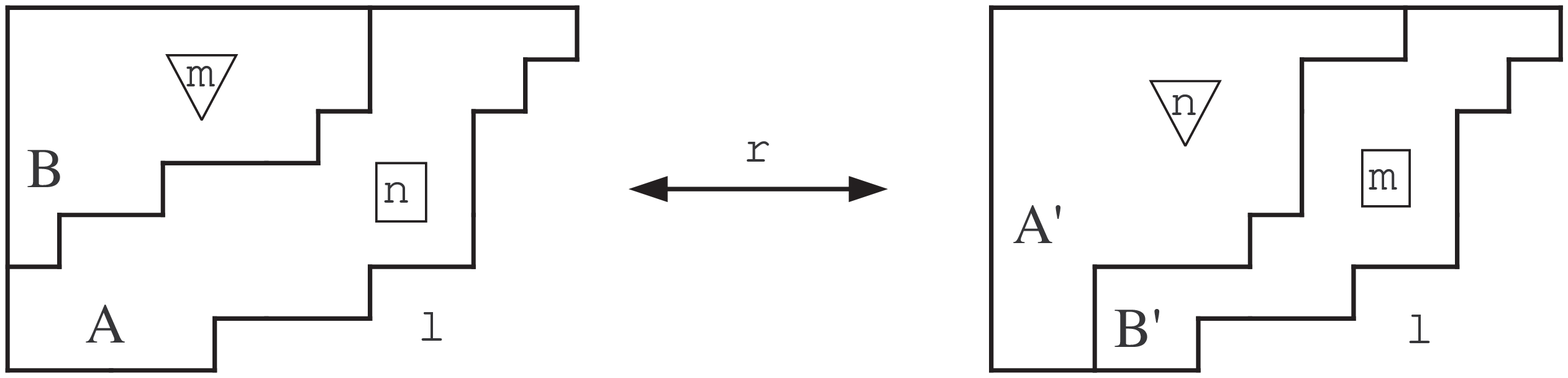,width=11cm}
\end{center}
\caption{Illustration of $\rho_1: A \to B^\pr$, where
$A \in \LR(\la/\mu,\nu)$, $B^\pr \in \LR(\la/\nu,\mu)$.}
\label{f:fundsym1}
\end{figure}

\begin{figure}[hbt]
\vskip.3cm
\psfrag{A}{$A$}
\psfrag{B}{$B$}
\psfrag{C}{$C$}
\psfrag{A'}{$A'$}
\psfrag{B'}{$B'$}
\psfrag{C'}{$C'$}
\psfrag{l}{$\la$}
\psfrag{m}{$\mu$}
\psfrag{n}{$\ba$}
\psfrag{s}{$\si$}
\psfrag{t}{$\tau$}
\psfrag{x}{$\zeta$}
\begin{center}
\epsfig{file=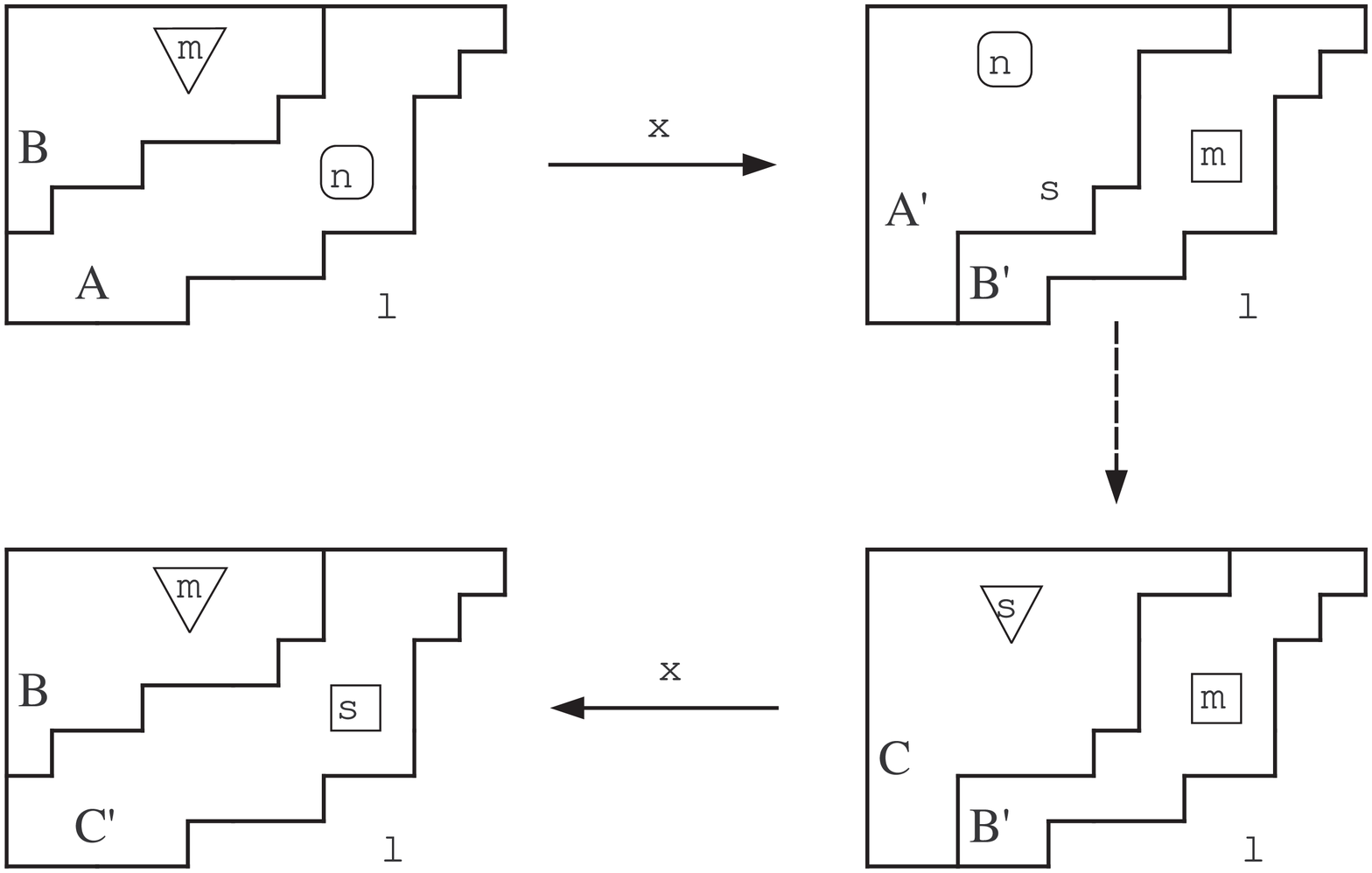,width=11cm}
\end{center}
\caption{Illustration of $\phi: A \to (A^\pr,C^\pr)$, where
$A \in \YT(\la/\mu,\ba)$, $A^\pr \in \YT(\si,\ba)$,
and $C^\pr \in \LR(\la/\mu,\si)$, for some $\si \vdash |\la/\mu|$.
}
\label{f:littlewood}
\end{figure}

We define the Littlewood-Robinson map~$\phi$ as follows:
Let $A\in\YT(\la/\mu, \ba)$ and $B=\Can(\mu)$.
Consider $(A^\prime, B^\prime)=\zeta(B,A)$, then
$A^\prime\in\YT(\sigma,\ba)$ for some $\sigma \vdash |\la/\mu|$.
Moreover, by Propositions~\ref{swijeu}, \ref{jeucan} and the
involution property~$(\maltese)$, $B^\prime\in\LR(\la/\si,\mu)$.
Now, let~$C = \Can(\si)$ and $(B^{\pr\pr},C^\pr) = \zeta(C,B^\pr)$.
Since $B^\pr \in\LR(\la/\si,\mu)$, we have
$B^{\pr\pr} = B = \Can(\mu)$ and therefore
$C^\pr \in \YT(\la/\mu,\si)$.
Similarly, since $C = \Can(\si)$,
we have~$C^\pr \in \LR(\la/\mu,\si)$.
Finally, let~$\phi(A) =  (A^\pr,C^\pr)$.
Observe that $(A^\pr,C^\pr)$ is in
$\bigcup_\pi \, \YT(\pi,\ba) \times
\LR(\la/\mu,\pi)$, as desired (see Figure~\ref{f:littlewood}).

%

We summarize:

\begin{prop}\label{liro}
The map~$\phi$ is a well defined bijection.
\end{prop}


\subsection{Second Fundamental Symmetry map}
\label{sec:prop-second}

We define $\rho_2 : \LR(\lambda/\mu,\nu) \to \LR(\lambda/\nu,\mu)$ as
the composition of two linear maps, $\gamma$ and $\tau$, with the
Sch\"utzenberger involution~$\xi$.
For this we need to introduce some notation.
Given partitions $\la$, $\mu$, $\nu$, such that $|\la|=|\mu| + |\nu|$,
we define
\[
\FV(\mu,\nu,\la) \, = \, \{ B\in\YT(\mu,\la -\nu) \mid
B \circ \Can(\nu) \in \LR(\mu \circ \nu,\la)\},
\]
\[
\FV^*(\mu,\nu,\la)\, = \, \{ B\in\YT(\mu, (\la -\nu)^*) \mid
B^\bu \circ \Can(\nu) \in \LR(\mu^\bu \circ \nu, \la )\}.
\]

\noindent
Now the map
\[
\gamma: \LR(\lambda/\mu,\nu) \to \FV^*(\mu, \nu, \la)
\]
is defined as follows.

Let $A \in \LR(\la/\mu,\nu)$ be a LR-tableau with recording matrix
$C=(c_{i,j})$, that is, $c_{i,j}$ is the number of $j$'s in the $i$-th row;
set also, $c_{0,0} = 0$, $c_{i,0} = \mu_i$ for
$1\le i\le l= \ell(\lambda)$ and $c_{k,j} = 0$ for $k>l$.
Define $A'=\gamma(A)$ to be the tableau with recording matrix
$D = (d_{i,j})$, where
\begin{figure}[hbt]
\vskip.3cm
\psfrag{A}{$A$}
\psfrag{B}{$A'$}
\psfrag{C}{$C = $}
\psfrag{D}{$D = $}
\psfrag{E}{$E = $}
\psfrag{F}{$F = $}
\psfrag{A'}{$A'$}
\psfrag{B'}{$B'$}
\psfrag{C'}{$C'$}
\psfrag{A''}{$A''$}
\psfrag{B''}{$B''$}
\psfrag{g}{$\ga$}
\psfrag{l}{$\la$}
\psfrag{m}{$\mu$}
\psfrag{n}{$\ba$}
\psfrag{h}{$\ba{\hspace{-.3ex}}^\ast$}
\psfrag{s}{$\si$}
\psfrag{x}{$\xi$}
%
\begin{center}
\epsfig{file=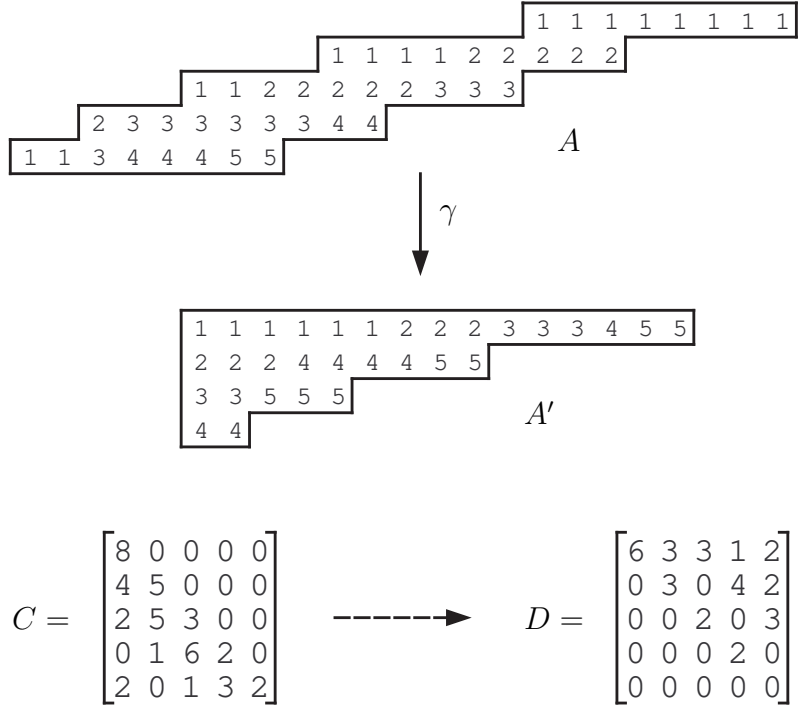,width=11cm}
\end{center}
\caption{Example of a LR-tableau $A \in \LR(\la/\mu,\nu)$,
the corresponding tableau $A' =\ga(A) \in \FV^\ast(\mu,\nu,\la)$,
and their recording matrices~$C$ and~$D$.
Here $\la = (23,18,15,11,8)$,
$\mu = (15,9,5,2,0)$, and $\nu = (16,11,10,5,2)$.}
\label{f:example-1}
\end{figure}
\[
d_{i,j} \, = \, \sum_{q=0}^{l-j} c_{l-j+i,\, q} \, -
\sum_{q=0}^{l-j+1} c_{l-j+1+i,\, q} \, .
\]
An example is given in Figure~\ref{f:example-1}.

\begin{prop} \label{carful}
The map $\gamma$ is a well defined bijection.
\end{prop}

The map
\[
\tau : \LR(\lambda/\mu, \nu) \to \FV(\nu, \mu, \la)
\]
is defined as follows:
Let $A \in \LR(\lambda/\mu,\nu)$ be a tableau with recording matrix\break
$C=(c_{i,j})$.
Define $\tau(A)$ to be the tableau with recording matrix
$E=(e_{i,j})$, where $e_{i,j} = c_{j,i}$, for $1\le i\le j\le l$.

\begin{prop}\label{tau}
The map $\tau$ is a well defined bijection.
\end{prop}

\begin{figure}[hbt]
\vskip.3cm
\psfrag{A}{$A$}
\psfrag{B}{$A'$}
\psfrag{C}{$C = $}
\psfrag{D}{$D = $}
\psfrag{E}{$E = $}
\psfrag{F}{$E' = $}
\psfrag{A'}{$A'$}
\psfrag{B'}{$B$}
\psfrag{C'}{$C'$}
\psfrag{Z}{$B'$}
\psfrag{B''}{$B''$}
\psfrag{C''}{$C''$}
\psfrag{l}{$\la$}
\psfrag{m}{$\mu$}
\psfrag{t}{$\tau^{-1}$}
\psfrag{x}{$\xi$}
%
\begin{center}

\epsfig{file=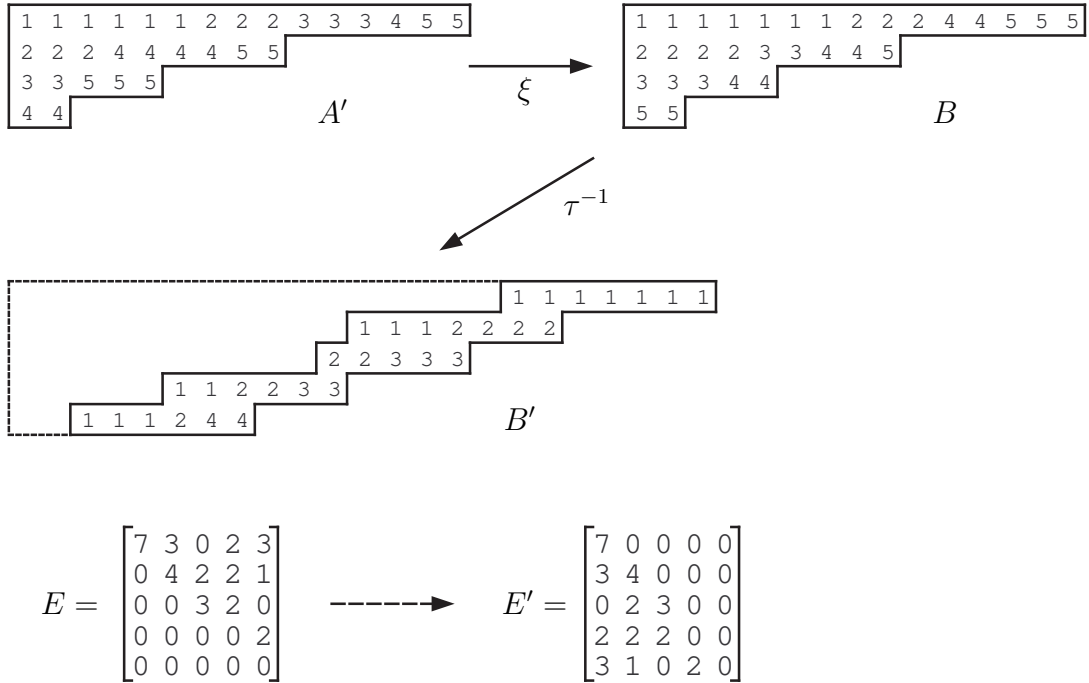,width=14.8cm}
\end{center}
\caption{Example of the Sch\"utzenberger involution $\xi: A' \mto B$
and the map $\tau^{-1}: B \mto B'$, and the recording matrices~$E, E'$
of the Young tableaux~$B,B'$.
Here $A' \in \FV^\ast(\mu,\nu,\la)$, $B \in \FV(\mu,\nu,\la)$,
and $B' \in \LR(\la/\nu,\mu)$, where $\la,\mu,\nu$
as in Figure~\ref{f:example-1}.}
\label{f:example-2}
\end{figure}

\bigskip
The map $\rho_2$ is defined as the composition
$\rho_2 := \tau^{-1} \circ \xi^\N \circ \gamma$.
A similar fundamental symmetry map can be defined as
the its inverse:
${\rho_2}^\prime := \rho_2^{-1} = \gamma^{-1} \circ \xi^\N \circ \tau$.


\begin{prop}\label{funsim2}
The maps $\rho_2$ and ${\rho_2}^\prime$ are
\emph{fundamental symmetry maps}.
\end{prop}

Note that, by Proposition~\ref{evajeu} and the involution property
$(\maltese)$, the restriction of the Sch\"utzenberger involution
defines a bijection
$\xi^\N : \FV^\ast(\mu,\nu,\la) \to \FV(\mu,\nu,\la)$.
This, together with Propositions~\ref{carful} and~\ref{tau} imply
Proposition~\ref{funsim2}.

\smallskip
We conjecture in Section~\ref{sec:further-alz} that $\rho_1$
in fact coincides with~$\rho_2$ and with ${\rho_2}^\pr$.
In the absence of this result
we use a different argument to prove linear equivalence
of these maps and the remaining maps in Theorem~1.

\bigskip
\subsection{Reversal}
\label{sec:rev}

\begin{figure}[hbt]
\vskip.3cm
\psfrag{A}{$A$}
\psfrag{B}{$B$}
\psfrag{C}{$C$}
\psfrag{D}{$\xi(A')$}
\psfrag{A'}{$A'$}
\psfrag{B'}{$B'$}
\psfrag{C'}{$C'$}
\psfrag{A''}{$A''$}
\psfrag{B''}{$B''$}
\psfrag{C''}{$C''$}
\psfrag{l}{$\la$}
\psfrag{m}{$\mu$}
\psfrag{n}{$\ba$}
\psfrag{h}{$\ba{\hspace{-.3ex}}^\ast$}
\psfrag{s}{$\si$}
\psfrag{x}{$\zeta$}
\psfrag{y}{$\xi$}
\begin{center}
\epsfig{file=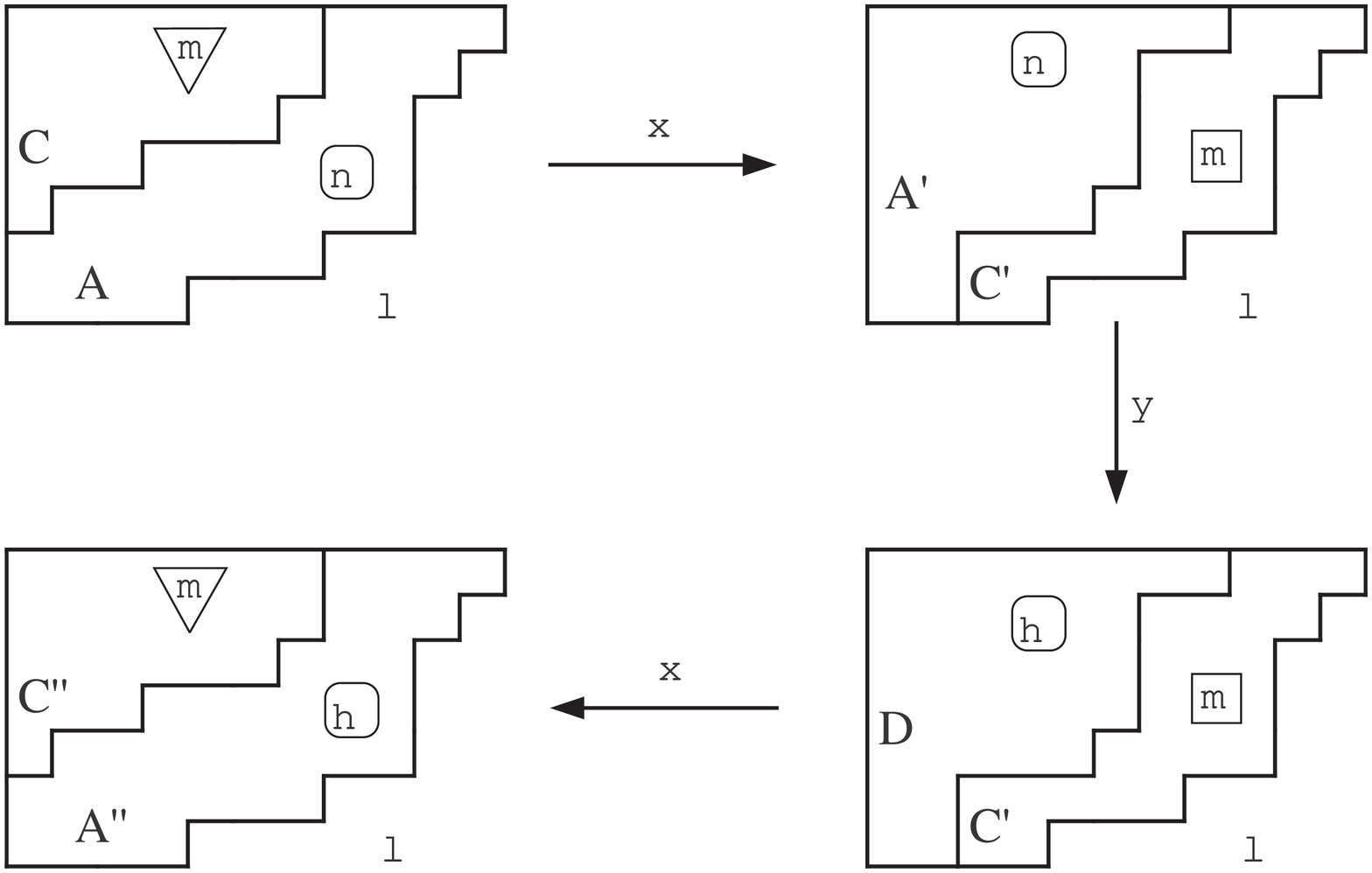,width=11cm}
\end{center}
\caption{Illustration of $\chi: A \mto A''$, where
$A \in \YT(\la/\mu,\ba)$, $A^{\pr\pr} \in \YT(\la/\mu,\ba^\ast)$.
}
\label{f:reversal}
\end{figure}

We define reversal as the {\em conjugation} of Sch\"utzenberger
involution with tableau switching $\chi=\zeta\circ\xi\circ\zeta$.
More precisely, let $A\in\YT(\la/\mu,\ba)$ and denote $C=\Can(\mu)$.
Consider $(A^\prime, C^\prime)=\zeta(C,A)$ and
$(C^{\prime\prime}, A^{\prime\prime})=\zeta(\xi(A^\prime),C^\prime)$.
Then, by Propositions~\ref{swijeu}, \ref{jeucan} and the
involution property~(\maltese), $C^{\prime\prime}= \Can(\mu)$,
and therefore $A^{\pr\pr}\in\YT(\la/\mu,\ba^\ast)$.
Reversal is defined by $\chi(A):=A^{\prime\prime}$
(see Figure~\ref{f:reversal}).

\begin{prop}
The map $\chi$ is a well defined involution.
\label{rev}
\end{prop}


\vskip1.cm

\section{Collection of linear reductions} \label{sec:zoo}

\subsection{Outline of the proof of Theorem~1}
\label{sec:zoo-outline}

Following definitions, we would have to present
$8\cdot 7 = 56$ different linear reductions to prove
Theorem~1.  In fact, as will be shown in the next section,
linear equivalence is an equivalence relation, so only
$2 \cdot (8-1) = 14$ linear reductions suffice
(say, between~$\vp$ and all other maps).
Some of these reductions are quite difficult, while
the proof can be obtained by a smaller number of
easier reductions.  The latter follows a less obvious
pattern summarized in the following lemma.

\smallskip

Recall that $\al \lrm \be$ stands for the map~$\al$ linearly
reducible to the map~$\be$.

\bigskip

{\bf Main Lemma} \ {\it For the maps as in Theorem~1, the
following linear reductions hold:

\smallskip

$\bu$ \ $\vp \lrm \psi \lrm \phi \lrm \zeta^\N \lrm \xi^\N \lrm \vp$,

\smallskip

$\bu$ \  $\rho_1 \lrm \zeta^\N \lrm \zeta \lrm \rho_1$ \, and \,
$\rho_2 \lrm \xi^\N \lrm \rho_2$,

\smallskip

$\bu$ \ $\chi \lrm \xi^\N \lrm \chi$.
}

\bigskip

In the next subsection we show that the Main Lemma implies
the first part of Theorem~1.
The rest of the section will contain the proof of linear
reductions in the Main Lemma.

\bigskip

\subsection{Compositions of linear reductions}
\label{sec:zoo-comp}

We start with the following elementary but very useful
results, which simplify the proof of Theorem~1.
While these results are standard in
Computer Science, they have never appeared in this
context.  We present short straightforward
proofs for completeness.

\bigskip

{\bf Composition Lemma} \ {\it
Suppose $\al_1 \lrm \al_2$ and $\al_2 \lrm \al_3$.
Then $\al_1 \lrm \al_3$. Moreover, if $\al_1$ can be
computed in at most~$s_1$ cost of~$\al_2$,
and $\al_2$ can be
computed in at most~$s_2$ cost of~$\al_3$, then
$\al_1$ can be
computed in at most~$(s_1s_2)$ cost of~$\al_3$.
}

\medskip

{\bf Proof.} \ Suppose $\CC_1$ is a $\al_2$-based ps-circuit
which computes~$\al_1$, and $\CC_2$ is a $\al_3$-based
ps-circuit which computes~$\al_2$. Substitute each of the
$\bs(\CC_1)$ maps~$\al_2$ in circuit~$\CC_1$ with the
circuit~$\CC_2$.  Denote the resulting circuit by~$\CC$.
By definition, $\CC$ is a $\al_3$-based ps-circuit
computing~$\al_1$. Observe also that~$\al_3$ is
used~$\bs(\CC_2)$ times in each copy of~$\CC_2$, and
thus~$\al_3$ is
used~$\bs(\CC_1)\bs(\CC_2)$ times in~$\CC$.  This
implies the result. \ $\sq$

\bigskip

{\bf Corollary 1.} \ {\it Suppose $\al_1 \sim \al_2$ and
$\al_2 \sim \al_3$.  Then $\al_1 \sim \al_3$.
}

\medskip

{\bf Proof.} \ By definition of linear equivalence, we have
$\al_1 \lrm \al_2 \lrm \al_3$.  Now the Composition Lemma
implies~$\al_1 \lrm \al_3$.  Similarly,
$\al_3 \lrm \al_2 \lrm \al_1$, which implies~$\al_3 \lrm \al_1$.
We conclude~$\al_1 \sim \al_3$.  \ $\sq$

\bigskip

{\bf Corollary 2. (Cycle Lemma)} \ {\it Suppose
$\al_1 \lrm \al_2 \lrm \ldots \lrm \al_n \lrm \al_1$.
Then $\al_1 \sim \al_2 \sim \ldots \sim \al_n$.
}

\medskip

{\bf Proof.} \ For every $1 \le i < j \le n$  we have
\, $\al_i \lrm \al_{i+1} \lrm \ldots \lrm \al_j$  \, and
$$\al_j \lrm \al_{j+1} \lrm \ldots \lrm \al_n \lrm \al_1
\lrm \al_2 \lrm \ldots \lrm \al_i
\,.$$
By Composition Lemma, this implies $\al_i \lrm \al_j$ and
$\al_j \lrm \al_i$, and thus~$\al_i \sim \al_j$. \ $\sq$

\bigskip

{\bf Corollary 3.} \ {\it Main Lemma implies the first part
of Theorem~1.
}

\medskip

{\bf Proof.} \ By Cycle Lemma, we have
$\vp \sim \psi \sim \phi \sim \zeta^\N \sim \xi^\N$.
These equivalences,
$\rho_1 \sim \zeta^\N \sim \zeta$, $\rho_2 \sim \xi^\N$,
$\chi \sim \xi^\N$ and Corollary~1 prove the claim.  \ $\sq$

\bigskip

\subsection{Reduction $\vp \lrm \psi$}
\label{sec:zoo-12}

The linear reduction of the RSK map~$\vp$ to the
Jeu de Taquin map~$\psi$ follows from Proposition~3.
Below we present the corresponding $\psi$-based
ps-circuit proving $\vp \lrm \psi$.

\medskip

$\di$ \ Input $k$, $\ba,\bb,$ such that~$\ell(\ba), \ell(\bb) \le k$.

$\di$ \ Input  $V = (v_{i,j}) \in \Mat(\ba,\bb)$.

$\di$ \ Set $\pi:= (a_1+ \cdots + a_k, a_1 + \cdots + a_{k-1},\ldots, a_1)$.

$\di$ \ Set $\sigma :=
( a_1 + \cdots + a_{k-1}, a_1 + \cdots + a_{k-2}, \dots, a_1, 0)$.

$\di$ \ Set $\rho:= (b_1+ \cdots + b_k, b_1 + \cdots + b_{k-1},\ldots, b_1)$.

$\di$ \ Set $\tau:=( b_1 + \cdots + b_{k-1}, b_1 + \cdots + b_{k-2}, \dots, b_1, 0)$.

$\di$ \ Set $V^{\updownarrow} :=( v_{k+1-i,j}) \in \Mat(\ba^*, \bb)$ and
$V^{\pr \updownarrow} := (v_{j,k+1-i}) \in \Mat(\bb^*,\ba)$.

$\di$ \ Let $Y \in \YT(\pi/\sigma,\bb)$ be the tableau with
recording matrix $V^\updownarrow$.

\qquad $\vdi$ \ Compute $B = \psi(Y)$.

$\di$ \ Let $X \in \YT(\rho/\tau,\ba)$ be the tableau with
recording matrix $V^{\pr \updownarrow}$.

\qquad $\vdi$ \ Compute $A = \psi(X)$.

$\di$ \ Output $(B,A) = \vp(V)$.

\medskip

\noindent
The above circuit is a simple parallel circuit which uses~$\psi$
twice. To prove that it computes~$\vp$, apply Proposition~\ref{rsk}.

\bigskip


\subsection{Reduction $\psi \lrm \phi$}
\label{sec:zoo-23}

The linear reduction of the Jeu de Taquin map~$\psi$
to the Littlewood-Robinson map~$\phi$ follows immediately from
Proposition~\ref{lirojeu}.  The corresponding circuit is a trivial
circuit~$I(\text{id},\phi,\de)$, where
$\text{id}$ is the identity map, and~$\de$
is a projection onto the first component.

\bigskip

\subsection{Reduction $\phi \lrm \zeta^\N$}
\label{sec:zoo-34}

The linear reduction of the Littlewood-Robinson map~$\phi$ to
the Tableau Switching map for normal shapes~$\zeta^\N$ follows from
the definition of $\phi$ given before Proposition~\ref{liro}.
Below we present the corresponding $\zeta^\N$-based ps-circuit
proving $\phi \lrm \zeta^\N$.

\medskip

$\di$ \ Input $k$, $\ba$, partitions~$\la, \mu$, such that
$\ell(\ba), \ell(\la)  \le k$.

$\di$ \ Input $A \in \YT(\la/\mu,\ba)$.

$\di$ \ Set $B = \Can(\mu)$.

\qquad $\vdi$ \ Compute $(A^\pr,B^\pr) = \zeta^\N(B,A)$.

$\di$ \ Let $\si$ be the shape of~$A^\pr$.

$\di$ \ Set $C = \Can(\si)$.

\qquad $\vdi$ \ Compute $(B^{\pr\pr},C^\pr) = \zeta^\N(C,B^\pr)$.

$\di$ \ Output $(A^\pr, C^\pr) \in \YT(\si,\ba)
\times \LR(\la/\mu,\si)$.

\medskip

\noindent
The above circuit is a simple sequential circuit which
uses map~$\zeta^\N$ twice.  Its correctness follows immediately
from our definition of $\phi$.

\bigskip

\subsection{Reductions $\zeta \lrm \xi$ and $\zeta^\N \lrm \xi^\N$}
\label{sec:zoo-45}


The linear reduction of the Tableau Switching
map~$\zeta$ to the Sch\"utzenberger involution~$\xi$
is given by the following simple sequential circuit.

\medskip

$\di$ \ Input $k$, $\mu \ssu \pi \ssu \la$, $\ba,\bb$, such that
$\ell(\la),\ell(\ba),\ell(\bb) \le k$.

$\di$ \ Input $A \in \YT(\la/\pi,\ba)$, $B \in \YT(\pi/\mu,\bb)$.


\qquad $\vdi$ \ Compute $B^\pr = \xi(B)\in
\YT(\pi/\mu, \bb^\ast)$, where $\bb^\ast = (b_k,\ldots,b_1)$.

$\di$ \ Relabel integers in~$A$ by adding~$k$ to them.

$\di$ \ Let $C := B^\pr \star A \in
\YT(\la/\mu, (b_k,\ldots, b_1, a_1,\ldots, a_k))$.

\qquad $\vdi$ \ Compute $C^\pr = \xi(C)\in
\YT(\la/\mu, (a_k,\ldots, a_1, b_1, \ldots, b_k))$.

$\di$ \ Decompose $C^\pr := A^{\pr} \star B^{\pr\pr}$, where
$B^{\pr\pr} \in \YT(\la/\si,\wt \bb )$,
$A^{\pr} \in \YT(\si/\mu,\ba^\ast)$,

\hskip .5cm
$\wt \bb  = (0,\ldots,0, b_1,\ldots, b_k)$, with $k$~zeros,
and $\si \vdash |\mu|+|\ba|$.

$\di$ \ Relabel integers in~$B^{\pr\pr}$ by subtracting~$k$ from them.

\qquad $\vdi$ \ Compute $A^{\pr\pr} = \xi(A^\pr)\in
\YT(\si/\mu, \ba)$.

$\di$ \ Output $(A^{\pr\pr},B^{\pr\pr}) = \zeta(B,A)$.

\medskip

\begin{figure}[hbt]
\vskip.3cm
\psfrag{A}{$A$}
\psfrag{B}{$B$}
\psfrag{C}{$C$}
\psfrag{A'}{$A'$}
\psfrag{B'}{$B'$}
\psfrag{C'}{$C'$}
\psfrag{A''}{$A''$}
\psfrag{B''}{$B''$}
\psfrag{C''}{$C''$}
\psfrag{l}{$\la$}
\psfrag{m}{$\mu$}
\psfrag{a}{$\ba$}
\psfrag{b}{$\bb$}
\psfrag{a^}{$\wt\ba$}
\psfrag{b^}{$\wt\bb$}
\psfrag{aa}{$\ba{\hspace{-.3ex}}^\ast$}
\psfrag{bb}{$\bb{\hspace{-.3ex}}^\ast$}
\psfrag{s}{$\si$}
\psfrag{p}{$\pi$}
\psfrag{x}{$\xi$}
\begin{center}
\epsfig{file=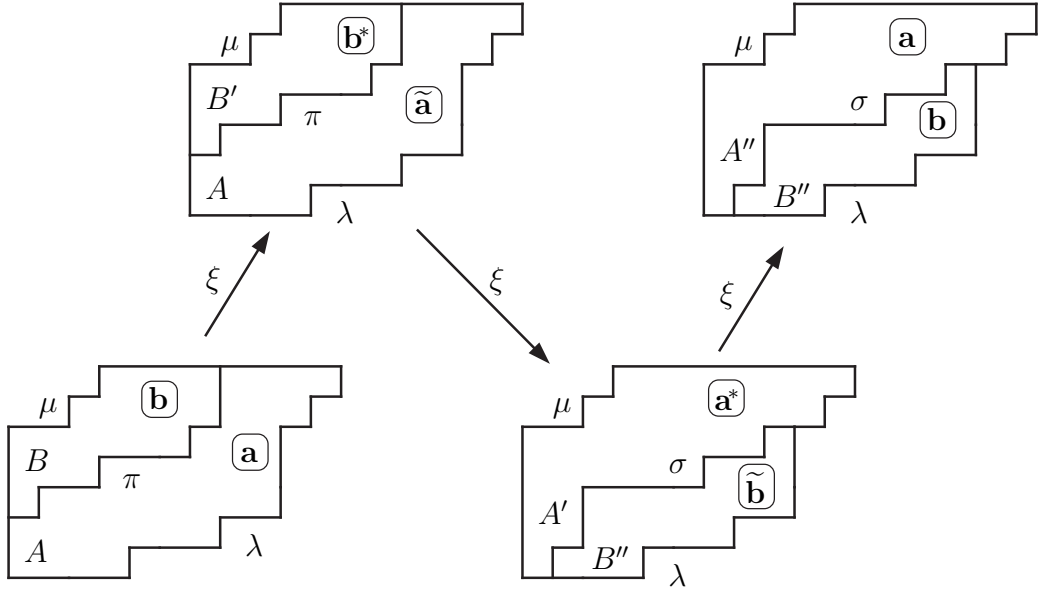,width=14.cm}
\end{center}
\caption{Illustration of the linear reduction $\zeta \lrm \xi$.
}
\label{f:involts}
\end{figure}

\noindent
This gives a sequential circuit which uses~$\xi$ three times.
We illustrate it in Figure~\ref{f:involts}.
Here we use $\wt\ba = (0,\ldots,0,a_1,\ldots,a_k)$,
with~$k$ zeros and keep notation~$A, B''$ for
tableaux before and after relabelling.  We hope this
won't lead to the confusion.


\medskip

\bigskip

{\bf Lemma 1.} \ {\it The above $\xi$-based circuit
computes~$\zeta$.}

\bigskip

\noindent
We postpone the proof till
Section~\ref{sec:proofs-lemma1}.
The proof is based on Propositions~\ref{evac} and~\ref{switch} in
Section~\ref{sec:prop-bk}.
When $\mu$ is the empty partition the above becomes a $\xi^\N$-based
circuit that computes $\zeta^\N$.

\bigskip

\subsection{Reduction $\xi^\N \lrm \vp$}
\label{sec:zoo-51}

The linear reduction of the Sch\"utzenberger
involution for normal shapes~$\xi^\N$ to the RSK correspondence~$\vp$
is given by the following simple sequential circuit.
The construction is based on Proposition~\ref{rskeva}.


\medskip

$\di$ \ Input $k$, $\la$, $\ba$, such that
$\ell(\la), \ell(\ba) \le k$.

$\di$ \ Input $A \in \YT(\la,\ba)$.




$\di$ \ Set $V = (v_{i,j})$ be the recording matrix of~$A$.
Let $U = (u_{i,j}):=(v_{i,k-j+1})$.


\qquad $\vdi$ \ Compute $(A^\circ,B^\circ) = \vp(U)$.



$\di$ \ Output $A^\circ = \xi^\N(A) \in \YT(\la,\ba^\ast)$.

\medskip

\noindent
The above is a trivial circuit which
uses map~$\vp$ only once.

\bigskip

{\bf Lemma 2.} \ {\it The above $\vp$-based circuit
computes~$\xi^\N$.}

\bigskip

The proof of the lemma follows from Proposition~\ref{rskeva} and additional
considerations, which will be given in Section~\ref{sec:proofs-lemma2}.

\bigskip

\subsection{Reductions $\rho_1 \lrm \zeta^\N \lrm \zeta$}
\label{sec:zoo-64}

The linear reduction of the First Fundamental Symmetry
map~$\rho_1$ to the Tableau Switching map for normal
shapes~$\zeta^\N$ follows immediately from Proposition~\ref{funsim1}.
The corresponding circuit
is a trivial circuit~$I(\de_1,\zeta,\de_2)$ where $\de_1$
creates a canonical tableau~$B=\Can(\mu)$,
$\zeta:(B,A) \mapsto (A^\pr,B^\pr)$,
and $\de_2$ is a projection on the second component~$B^\pr$.
We leave the easy technical details to the reader.
The linear reduction $\zeta^\N \lrm \zeta$ is trivial.

\bigskip

\subsection{Reduction $\zeta  \lrm \rho_1$}
\label{sec:zoo-46}

This reduction is more involved than other
linear reductions, and requires
an intermediate map~$\zeta^\LRN$. Formally, we first present
a linear reduction $\zeta \lrm \zeta^\LRN$, and then
a linear reduction $\zeta^\LRN  \lrm \rho_1$.
Now the Composition Lemma gives the desired construction.

\subsubsection{Reduction $\zeta \lrm \zeta^\LRN$ \label{sec:zoo-461}}

Suppose $\mu \ssu \la$, $n = |\la/\mu|$, and $|\nu| + |\tau|=n$.
Define \emph{LR-Tableau Switching} to be a one-to-one correspondence:
$$\zeta^\LRN: \, \bigcup_{\pi \, \vdash |\la|-|\nu|} \,
\LR(\pi/\mu,\tau) \times \LR(\la/\pi,\nu) \ \too
\bigcup_{\si \, \vdash |\la|-|\tau|} \,
\LR(\si/\mu,\nu) \times \LR(\la/\si,\tau) ,$$
which is given by restriction of~$\zeta$ to the sets as above.

\begin{prop}\label{swilr}
The map~$\zeta^\LRN$ is a well defined bijection.
\end{prop}


The reduction~$\zeta^\LRN \lrm \zeta$ is trivial, but will not
be needed.  Below we show that $\zeta \lrm \zeta^\LRN$,
which implies that $\zeta^\LRN \sim \zeta$.

\medskip


We first describe the working of the reduction.
Start with tableaux~$A$, $B$ and consider a tableau
$(B\star A) \comp C \comp D$, where $C$ contains only
integers as in $A$ and $D$ contains only
integers as in $B$ (see Figure~10).
Let $\wh A := A \comp_{s,0} C$, $\wh B:= B \comp_{t,k} D$ be
parts of the tableau above, for some $s$, $t$, $k$ to be defined;
recall the definition of $\circ_{a,b}$ in Section~\ref{sec:basic}.
Clearly, tableau switching of $\wh A$ and $\wh B$ gives
$B^\pr \comp_{t,k} D$, $A^\pr \comp_{s,0} C$,
where $(A^\pr,B^\pr) = \zeta(B,A)$.
Now, if $A \comp_{s,0} C$ and $B \comp_{t,k} D$ are $\LR$-tableaux,
this gives a linear reduction as desired.
Below we show that one always find tableaux $C$, $D$ as above.

\begin{figure}[hbt]
\vskip.3cm
\psfrag{A}{$A$}
\psfrag{B}{$B$}
\psfrag{C}{$C$}
\psfrag{D}{$D$}
\psfrag{A'}{$A'$}
\psfrag{B'}{$B'$}
\psfrag{C'}{$C'$}
\psfrag{A''}{$A''$}
\psfrag{B''}{$B''$}
\psfrag{C''}{$C''$}
\psfrag{l}{$\la$}
\psfrag{m}{$\mu$}
\psfrag{a}{$\al$}
\psfrag{b}{$\be$}
\psfrag{c}{$\mathbf c$}
\psfrag{d}{$\mathbf d$}
\psfrag{s}{$\si$}
\psfrag{p}{$\pi$}
\psfrag{z}{$\zeta^\LRN$}
%
\begin{center}
\epsfig{file=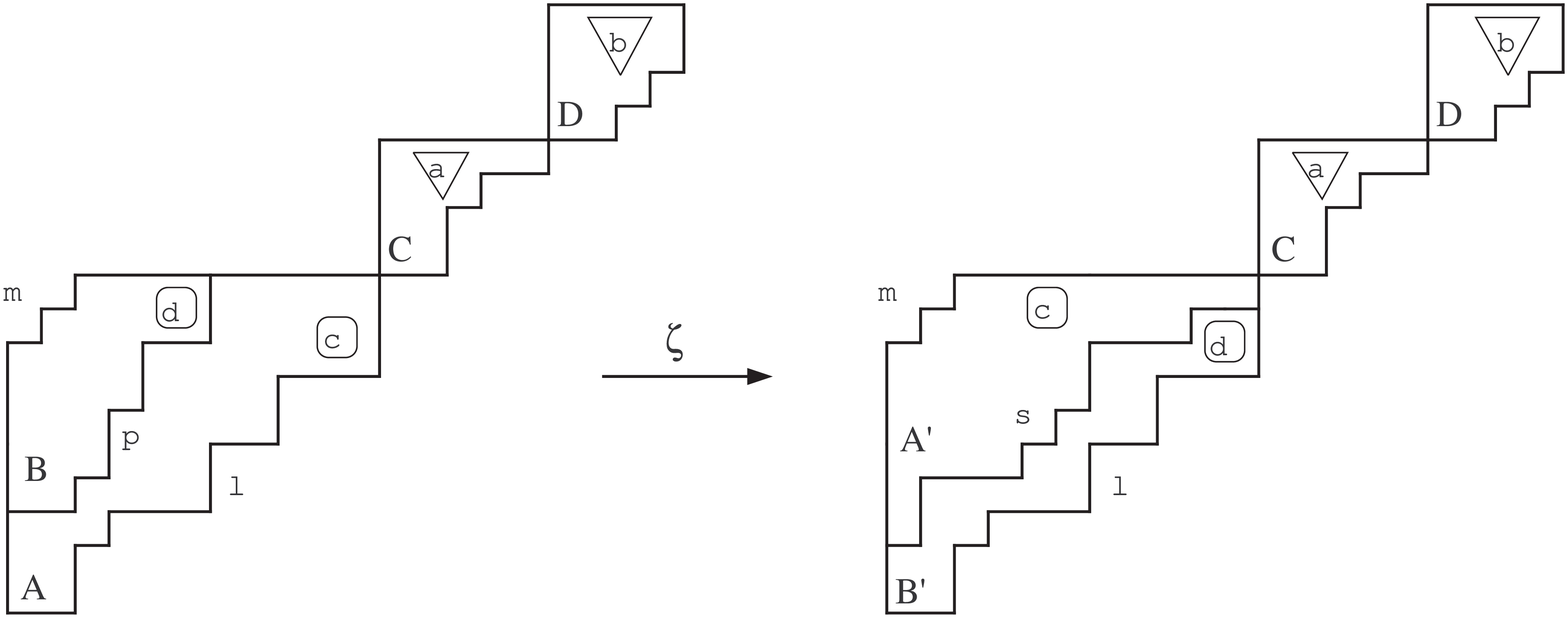,width=15.2cm}
\end{center}\label{f:lrts}
\caption{Illustration of the linear reduction $\zeta \lrm \zeta^\LRN$.
}
\end{figure}

\medskip

$\di$ \ Input $k$,  $\mu \ssu \pi \ssu \la$,
${\mathbf c}$, ${\mathbf d}$, such that
$\ell(\la),\ell({\mathbf c}),\ell({\mathbf d}) \le k$.

$\di$ \ Input $A \in \YT(\la/\pi, {\mathbf c})$,
$B \in \YT(\pi/\mu, {\mathbf d})$.

$\di$ \ Set $\al := (c_2 +\cdots + c_k,\, c_3+ \cdots + c_k,\ldots,\, c_k,\, 0)$.

$\di$ \ Set $\be := (d_2 +\cdots + d_k,\, d_3+ \cdots + d_k,\ldots,\, d_k,\, 0)$.


$\di$ \ Set $\wh \la:= (\la_1+\al_1+\be_1, \dots,\la_1+\al_1+\be_k,\,
\la_1+\al_1, \dots, \la_1+\al_k,\, \la_1, \dots, \la_k)$.

$\di$ \ Set $\wh \pi:= (\la_1+\al_1+\be_1, \dots,\la_1+\al_1+\be_k,\,
\la_1, \dots, \la_1,\, \pi_1, \dots, \pi_k)$.

$\di$ \ Set $\wh \mu:= (\la_1+\al_1, \dots,\la_1+\al_1,\,
\la_1, \dots, \la_1,\, \mu_1, \dots, \mu_k)$.

$\di$ \ Set $\wh A := A \comp_{\la_1,0} \Can(\al)
\in \LR(\wh \la/\wh \pi)$, \,
$\wh B := B \comp_{\la_1 +\al_1, k} \Can(\be) \in \LR(\wh \pi/\wh \mu)$.

\qquad $\vdi$ \ Compute $({\wh A}^\pr,{\wh B}^\pr) = \zeta^\LRN(\wh B,\wh A)$.

$\di$ \ Decompose ${\wh A}^\pr = A^\pr \comp_{\la_1,0} \Can(\al)$,
\, ${\wh B}^\pr = B^\pr \comp_{\la_1 + \al_1,k} \Can(\be)$, where

\hskip .5cm
$A^\pr \in \YT(\si/\mu, {\mathbf c})$, $B^\pr \in \YT(\la/\si, {\mathbf d})$,
for some $\si \vdash |\la/\pi|$.

$\di$ \ Output $(A^{\pr},B^{\pr}) = \zeta(B,A)$.

\bigskip

{\bf Lemma 3.} \ {\it The above} $\zeta^\LRN$-{\it based trivial circuit
computes~$\zeta$.}

\bigskip

Its correctness follows from Proposition~\ref{swilr} and
the fact that the Tableaux Switching map commutes with taking
compositions~\cite[Thm.~2.3]{BSS}.
We postpone the proof till Section~\ref{sec:proofs-lemma3}.

\subsubsection{Reduction $\zeta^\LRN \lrm \rho_1$ }
\label{sec:zoo-462}


First, we describe the working of the reduction.
Start with tableaux~$A \in \LR(\la/\pi,\nu)$, $B\in \LR(\pi/\mu,\tau)$.
We need to obtain $(A^\pr, B^\pr) = \zeta(B,A)$.
Think of~$\rho_1$ as tableau switching with a canonical tableau.
Computing $\rho_1(B)$
is a tableau switching of $\Can(\mu)$ and~$B$.
Similarly, computing
$\rho_1(\rho_1(B)\star A)$ is a tableau switching
which first returns to~$B$ and then switches~$B$ and~$A$
(see Figure~11).
This step works only if $\rho_1(B)\star A$ is a LR-tableau.
Thus, we extend $\rho_1(B)\star A$ to a tableau $E$ that is a LR-tableau,
and apply $\rho_1$ to $E$ instead of applying it to $\rho_1(B)\star A$.
Finally, we apply $\rho_1$ to restriction of tableau~$\rho_1(E)$
to integers~$1,\dots, k$ (see Figure~12).

\medskip

$\di$ \ Input $k$,  $\mu \ssu \pi \ssu \la$, $\nu,\tau$, such that
$\ell(\la),\ell(\nu),\ell(\tau) \le k$.

$\di$ \ Input $A \in \LR(\la/\pi,\nu)$, \, $B \in \LR(\pi/\mu,\tau)$.

\qquad $\vdi$ \ Compute $C:=\rho_1(B) \in \LR(\pi/\tau,\mu)$.

$\di$ \ Set $s:=\nu_1$, $t:=\la_1$.

$\di$ \ Set $\ga := (s^k)$, \, $G := \Can(\ga)$ .

$\di$ \ Set $\wh \pi := (t+s,\ldots,t+s,\pi_1,\ldots,\pi_k)$, \,
$\wh \la := (t+s,\ldots,t+s,\la_1,\ldots,\la_k)$,

\hskip .5cm
$\wh \tau := (t,\ldots,t,\tau_1,\ldots,\tau_k)$,
all of length $2k$.

$\di$ \ Set $\wt \mu := (s+\mu_1,\ldots,s+\mu_k)$, \,
$\vk := (s+\mu_1,\ldots,s+\mu_k,\, \nu_1,\ldots,\nu_k)$.

$\di$ \ Relabel the integers in~$A$ by adding~$k$ to them.

$\di$ \ Set $\wh C:= C \comp_{\la_1,0} G \in \LR(\wh \pi/\wh \tau,\wt \mu)$, \,
$D:=  \wh C \star A \in \LR(\wh \la/\wh \tau,\vk)$.

\qquad $\vdi$ \ Compute
$E:=\rho_1(D) \in \LR(\wh \la/\vk, \wh \tau)$.

$\di$ \ Set $\de:=(t,\dots, t)$ of length $k$, and
$\wt \tau = (0,\ldots,0,\tau_1,\ldots,\tau_k)$ of length $2k$.

$\di$ \ Decompose $E = F \star B^\pr$, where
$F\in\LR(\wh{\si}/\vk,\de)$,
$B^\pr\in\YT(\la/\si,\wt \tau)$,

\hskip .5cm
$\wh{\si} =(t+s,\dots, t+s,\, \si_1,\dots, \si_k)$ of
length~$2k$, for some $\si \vdash |\la/\tau|$.


$\di$ \ Relabel the integers in~$B^\pr$ by subtracting~$k$ from them:
now~$B^\pr\in\LR(\la/\si,\tau)$.

\qquad $\vdi$ \ Compute $H:=\rho_1(F)
\in\LR(\wh{\si}/\de,\vk)$.

$\di$ \ Set $\wt \nu = (0,\ldots,0,\nu_1,\ldots,\nu_k)$ of length $2k$.

$\di$ \ Decompose
$H = (\Can(\mu) \circ_{\la_1,0} G) \star A^\pr$,
where $A^\pr \in \YT(\si/\mu,\wt \nu)$.


$\di$ \ Relabel the integers in $A^\pr$ by subtracting $k$ from
them: now $A^\pr \in \LR(\si/\mu,\nu)$.

$\di$ \ Output $(A^\pr, B^\pr) \in \LR(\si/\mu,\nu) \times
\LR(\la/\si,\tau)$.

\medskip

\noindent
The above circuit is a sequential circuit which
uses map~$\rho_1$ three times.
Its correctness is summarized in the following lemma.

\bigskip

{\bf Lemma 4.} \ {\it The above $\rho_1$-based sequential circuit
computes the restricted tableaux switching map~{\rm $\zeta^\LRN$}.}

\bigskip

We prove the lemma in Section~\ref{sec:proofs-lemma4}.

\medskip

\begin{figure}[hbt]
\vskip.3cm
\psfrag{A}{$A$}
\psfrag{B}{$B$}
\psfrag{C}{$C$}
\psfrag{D}{$D$}
\psfrag{G}{$G$}
\psfrag{A'}{$A'$}
\psfrag{B'}{$B'$}
\psfrag{Ah}{$D$}
\psfrag{Bh}{$E$}
\psfrag{At}{$\wt A$}
\psfrag{Bt}{$F$}
\psfrag{C'}{$C'$}
\psfrag{a}{$\al$}
\psfrag{b}{$\be$}
\psfrag{det}{$\de=(t^k)$}
\psfrag{de}{$\de$}
\psfrag{k}{$k$}
\psfrag{ka}{$\vk$}
\psfrag{l}{$\la$}
\psfrag{m}{$\mu$}
\psfrag{n}{$\nu$}
\psfrag{g}{$\ga$}
\psfrag{p}{$\pi$}
\psfrag{s}{$\si$}
\psfrag{t}{$\tau$}
\psfrag{tt}{$\wh\tau$}
\psfrag{ttt}{$t$}
\psfrag{r}{$\rho_1$}
\psfrag{ss}{$s$}
\begin{center}
\epsfig{file=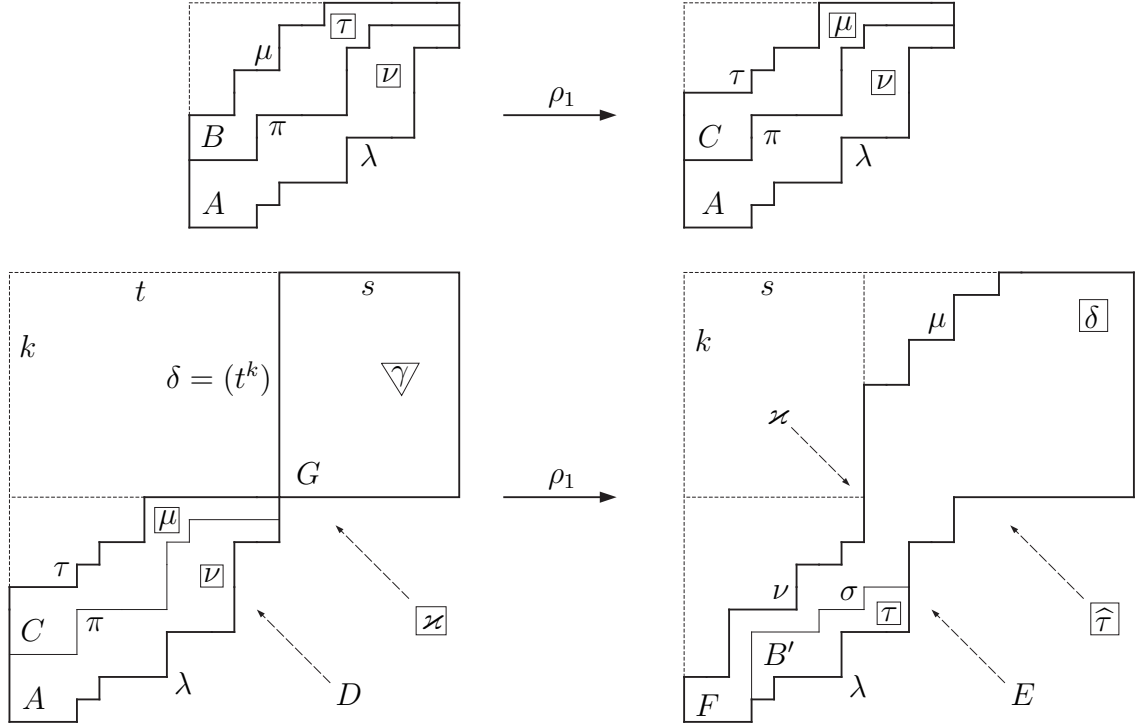,width=15.2cm}
\end{center}
\caption{Illustration of the first two applications of $\rho_1$
in linear reduction~$\zeta^\LRN \lrm \rho_1$.}
\label{f:tsts}
\end{figure}

\begin{figure}[hbt]
\vskip.3cm
\psfrag{A}{$A$}
\psfrag{B}{$B$}
\psfrag{C}{$C$}
\psfrag{D}{$D$}
\psfrag{G}{$G$}
\psfrag{A'}{$A'$}
\psfrag{B'}{$B'$}
\psfrag{Ah}{$\wh A'$}
\psfrag{Bh}{$\wh B'$}
\psfrag{At}{$H$}
\psfrag{Bt}{$F$}
\psfrag{C'}{$C'$}
\psfrag{a}{$\al$}
\psfrag{b}{$\be$}
\psfrag{det}{$\de=(t^k)$}
\psfrag{de}{$\de$}
\psfrag{k}{$k$}
\psfrag{ka}{$\vk$}
\psfrag{l}{$\la$}
\psfrag{m}{$\mu$}
\psfrag{n}{$\nu$}
\psfrag{nt}{$\wt\nu$}
\psfrag{g}{$\ga$}
\psfrag{p}{$\pi$}
\psfrag{s}{$\si$}
\psfrag{t}{$\tau$}
\psfrag{tt}{$\wh\tau$}
\psfrag{ttt}{$t$}
\psfrag{r}{$\rho_1$}
\psfrag{ss}{$s$}
\smallskip
\begin{center}
\epsfig{file=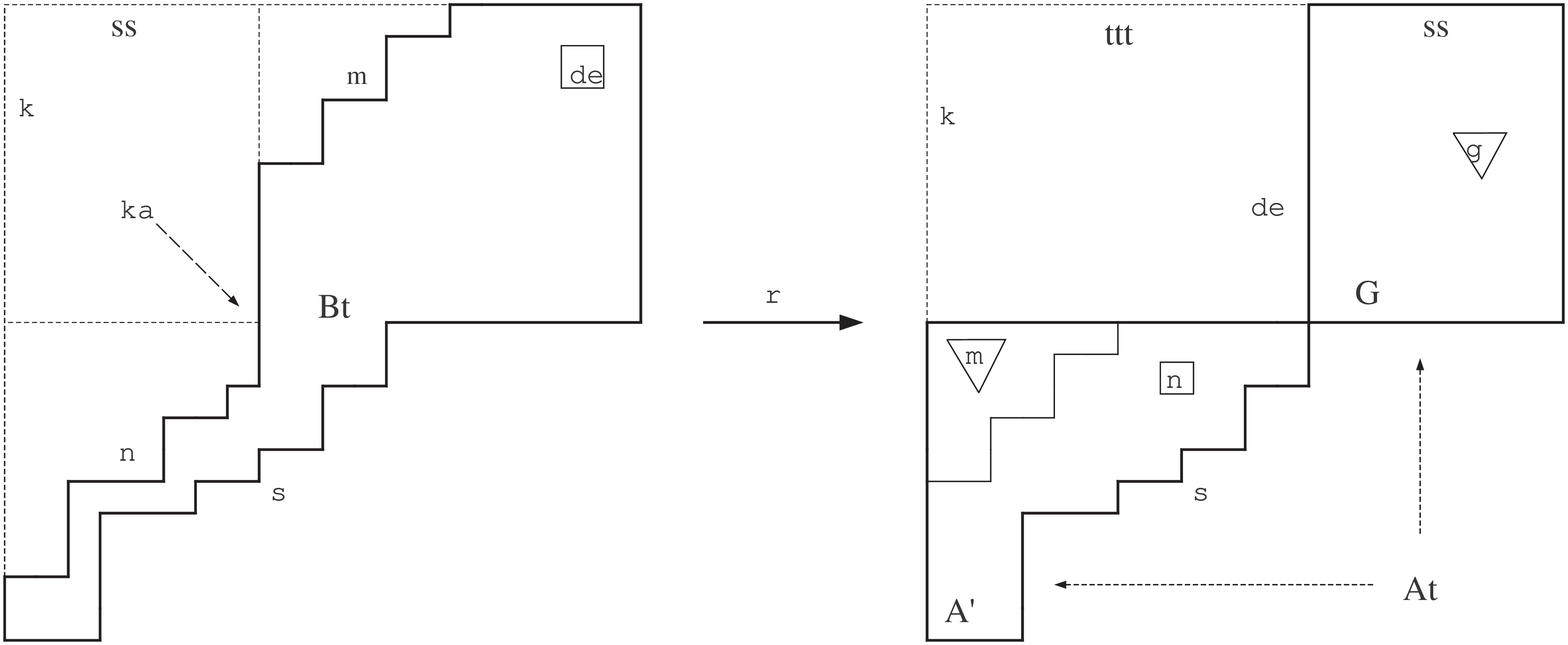,width=15.2cm}
\end{center}
\caption{Illustration of the of the third application
of~$\rho_1$ in linear reduction~{\rm $\zeta^\LRN \lrm \rho_1$}.
}
\label{f:tststs}
\end{figure}


\bigskip

\subsubsection{Using duality}
\label{sec:zoo-464}  One can use the duality
(rotating the picture 180 degree and relabelling
the integers) and apply~$\rho_1$ twice as in the beginning
of the circuit above (in place of the the third application
of~$\rho_1$).  This gives a conceptually easier $\rho_1$-based
ps-circuit for $\zeta^\LRN$, but with a higher cost.

\bigskip

\subsubsection{Reduction $\zeta \lrm \zeta^\N$ }
\label{sec:zoo-463}

Note that it is not difficult to define directly a linear
reduction $\zeta \lrm \zeta^\N$.  Even though we do not
need this reduction, let us quickly outline it.

Let $B\in\YT(\pi/\mu, {\mathbf d})$,
$A\in\YT(\la/\pi, {\mathbf c})$ and $C=\Can(\mu)$.
Relabel the entries of $B$ by adding $k$ to them, and compute
$(A^\pr, D^\pr)=\zeta^\N(C\star B, A)$.
Decompose $D^\pr=C^\pr \star B^\pr$, where $C^\pr$ has content
$\mu$.
Let $\sigma$ be the shape of $D^\pr$.
Relabel the entries of $B^\pr$ by subtracting $k$ from them;
thus $B^\pr\in\YT(\la/\si, {\mathbf d})$.
Let $(C^{\pr\pr}, A^{\pr\pr}) = \zeta^\N(A^\pr, C^\pr)$.
Since $C^{\pr\pr} =\Can(\mu)$, we have that
$A^{\pr\pr}\in\YT(\si/\mu, {\mathbf c})$.
We leave as an exercise to the interested reader to show that
$(A^{\pr\pr}, B^\pr)=\zeta(B,A)$.
In this way we obtain a sequential $\zeta^\N$-based
circuit which computes $\zeta$.


\bigskip

\subsection{Reduction $\rho_2 \lrm \xi^\N$}
\label{sec:zoo-75}

The linear reduction of the Second Fundamental Symmetry
map~$\rho_2$ to Sch\"utzenberger involution map~$\xi$
follows immediately from the definition of $\rho_2$
and Proposition~14.
The corresponding circuit is a trivial
circuit~$I(\gamma,\xi,\tau^{-1})$.

\bigskip

\subsection{Reduction $\xi^\N  \lrm \rho_2$}
\label{sec:zoo-57}


We describe first the working of the reduction.
Let $A\in \YT(\mu,\ba)$, with $\ba=(a_1, \dots, a_k)$.
Define $\nu =(a_1+\cdots + a_{k-1}, a_1+\cdots +a_{k-2}, \dots, a_1,0)$
and $\la =\nu+\ba^*$.
Then $A^\bu\circ\Can(\nu)\in \LR(\mu^\bu\circ\nu,\la)$
and $A\in\FV^*(\mu,\nu,\la)$.
Let $ A^\circ = \tau\circ\rho_2\circ\gamma^{-1}(A)\in\FV(\mu,\nu,\la)$.
By the definition of~$\rho_2$ we have~$A^\circ=\xi^\N(A)$.
Therefore, we have proved that the following trivial circuit
$I(\gamma^{-1}, \rho_2,\tau)$ computes $\xi^\N$.

\medskip

$\di$ \ Input $k$, $\ba$, $\mu$, such that $\ell(\mu)$, $\ell(\ba)\le k$.

$\di$ \ Input $A\in \YT(\mu,\ba)$.

$\di$ \ Set $\nu:=(a_1+\cdots + a_{k-1}, a_1+\cdots +a_{k-2}, \dots, a_1,0)$,
$\la:=\nu+\ba^*$.

$\di$ \ Set $B:=\gamma^{-1} (A)$.

\qquad $\vdi$ \ Compute $C:=\rho_2 (B)$.

$\di$ \ Set $A^\circ := \tau (C)$.

$\di$ \ Output $A^\circ = \xi^\N(A) \in \YT(\mu, \ba^*)$.


\bigskip
\subsection{Reduction $\xi^\N \lrm \chi$}
\label{sec:zoo-58}

The reduction $\xi^\N \lrm \chi$ is trivial, since when $\mu = \varnothing$,
$\chi$ coincides with $\xi^\N$.

\bigskip
\subsection{Reduction $\chi \lrm \xi^\N$}
\label{sec:zoo-85}


The linear reduction of the reversal map to the Sch\"utzenberger
involution for normal shapes is given by the following
circuit.

\medskip
$\di$ \ Input $k$, $\mu \ssu \la$, $\ba$, such that $\ell(\la)$,
$\ell(\ba) \le k$.

$\di$ \ Input $A\in\YT(\la/\mu, \ba)$.

$\di$ \ Set $C:= \Can(\mu)$.

\qquad $\vdi$ \ Compute $C^\circ :=\xi^\N(C)  \in \YT(\mu, \mu^\ast)$.

$\di$ \ Set ${\mathbf c}:=(\mu_k,\dots,\mu_1,\, a_1,\dots, a_k)$,
${\mathbf d} := (a_1,\dots,a_k,\, \mu_1,\dots,\mu_k)$
of length $2k$.

$\di$ \ Set
$\tilde\mu :=(0,\dots, 0,\, \mu_1,\dots,\mu_k)$ of length $2k$.

$\di$ \ Relabel de entries of $A$ by adding $k$ to them.

$\di$ \ Let $B:= C^\circ \star A\in \YT(\la, {\mathbf c})$.

\qquad $\vdi$ \ Compute $B^\circ :=
\xi^\N(B)\in \YT(\la, {\mathbf c}^\ast)$

$\di$ \ Decompose $B^\circ = A^\circ \star C^\pr$, where
$A^\circ\in \YT(\nu, \ba^\ast)$,
$C^\pr \in \YT(\la/\nu, \tilde\mu)$ for some

\hskip .5cm
partition $\nu\vdash |\ba|$.

\qquad $\vdi$ \ Compute $A^{\circ\circ}:= \xi^\N(A^\circ)\in \YT(\nu,\ba)$.

$\di$ \ Let $B^{\circ\circ} := A^{\circ\circ}\star C^\pr \in
\YT(\la, {\mathbf d})$.

\qquad $\vdi$ \ Compute $B^{\bw} := \xi^\N(B^{\circ\circ})
\in \YT(\la, {\mathbf d}^\ast)$.

$\di$ \ Decompose $B^{\bw}:= C^\circ\star A^{\bw}$.

$\di$ \ Relabel the entries of $A^{\bw}$ by subtracting $k$
from them.

$\di$ \ Output $\chi(A):= A^{\bw}\in \YT(\la/\mu, \ba^\ast)$.

\begin{figure}[hbt]
\vskip.3cm
\psfrag{A}{$A$}
\psfrag{B}{$B$}
\psfrag{C}{$C$}
\psfrag{D}{$\xi(A')$}
\psfrag{A'}{$A'$}
\psfrag{B'}{$B'$}
\psfrag{C'}{$C'$}
\psfrag{Co}{$C^\circ$}
\psfrag{Ao}{$A^\circ$}
\psfrag{Aoo}{$A^{\circ\circ}$}
\psfrag{A''}{$A^\bw$}
\psfrag{l}{$\la$}
\psfrag{m}{$\mu$}
\psfrag{nu}{$\nu$}
\psfrag{m'}{$\mu{\hspace{-.2ex}}^\ast$}
\psfrag{n}{$\ba$}
\psfrag{h}{$\ba{\hspace{-.3ex}}^\ast$}
\psfrag{s}{$\si$}
\psfrag{x}{$\xi^\N$}
%
\begin{center}
\epsfig{file=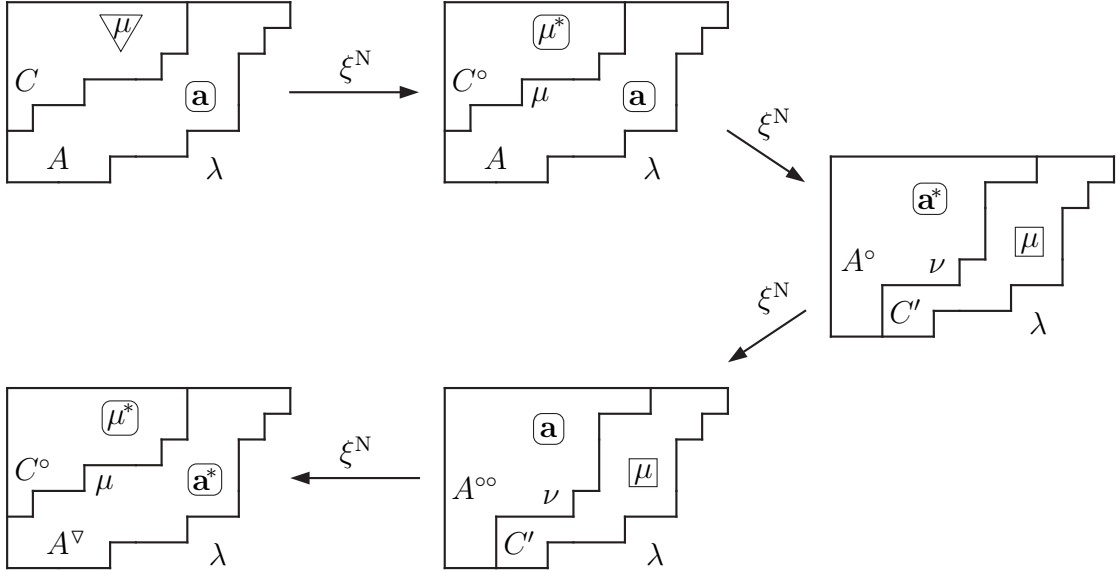,width=15cm}
\end{center}
\caption{Illustration of the reduction $\chi \lrm \xi^\N$.
}
\label{f:revxi}
\end{figure}

\medskip
The above circuit is a simple sequential circuit which uses the map
$\xi^\N$ four times (see Figure~\ref{f:revxi}).

\medskip
{\bf Lemma 5.} \ {\em The above $\xi^\N$-based circuit computes $\chi$.}

\medskip

In fact, since $C^\circ=\xi^\N(\Can(\mu))$ can be given explicitly,
we could write the circuit using the map~$\xi^\N$ only three times.
Formally, to prove this we need the following result:

\medskip
{\bf Lemma 6.} \ {\em Map $\mu \to \xi^\N(\Can(\mu))$ can be computed
in $O(k^2)$ time, where $k = \ell(\mu)$.}

\medskip

We prove both lemmas in Section~\ref{sec:proofs-lemmas}.

\bigskip

\subsection{Proof of Main Lemma}
\label{sec:zoo-together}

Main Lemma (Section~\ref{sec:zoo-outline})
now follows immediately from
reductions above (sections~\ref{sec:zoo-12} to~\ref{sec:zoo-85}).

\bigskip




\section{Proof of results}
\label{sec:proofs}

As we mentioned in the Introduction,
there is an extensive literature in
the subject of Young tableau bijections.  Thus in many
cases the technical results we need are already known.
In the next subsection we give a brief overview the
literature giving pointers to propositions.
Readers interested in a modern treatment of the subject,
complete definitions and further references, are referred
to~\cite{Fu,Sa,St}.  An historic overview of several of
the constructions used here can be found in~\cite{L1,L2}.
The remainder of this section contains proofs of lemmas
and theorems.

\subsection{Brief overview of the literature}
\label{sec:proofs-prop}

The Sch\"utzenberger involution or evacuation and the Jeu de Taquin
were introduced by\break \mbox{M.-P.~Sch\"ut}\-zen\-berger in his
study of the combinatorics of Young tableaux~\cite{Sc1,Sc3}.
The Sch\"ut\-zen\-berger involution
is usually considered only in the case $\mu=\emp$, but in fact
Sch\"ut\-zen\-berger extended this map to every poset~\cite{Sc}.
Proposition~\ref{evajeu} follows from Proposition~\ref{rsk}
and the duality theorem in Appendix A1 in \cite{Fu} or from
Theorem~A1.2.10 in~\cite{St}.
For Proposition~\ref{jeucan} see \cite[\S 5.2]{Fu} or \cite[A1.3.6]{St}.
An alternative extension of the Sch\"utzenberger involution to
skew shapes is the reversal map considered by Benkart, Sottile and
Stroomer in~\cite[\S 5]{BSS}.
The reversal of a skew tableau can be characterized by means of
{\em dual equivalence}~\cite{Ha},
although this is not the case for Sch\"utzenberger
involution in the case of skew shapes.
Proposition~\ref{rev} is proved in \S 5 of~\cite{BSS} and
follows from the fact that it is defined as the composition of
three bijections and relations $(\maltese)$.

Bender-Knuth transformations were introduced by E.~Bender and
D.E.~Knuth in their work on the enumeration of plane partitions
\cite{BK}.
Their connection with Sch\"utzenberger involution
was realized much later.
Proposition~\ref{evac} is explicit in~\cite[\S 2C]{KB}, although
it was probably discovered much earlier.

The RSK correspondence in full generality is due to
Knuth~\cite{Kn}, who built it based on previous works of
Robinson~\cite{Ro} and Schensted~\cite{Sch}
(and conversations with Sch\"utzenberger).
Proposition~\ref{rsk} can be found in \cite{Sc3} or \cite[\S 4]{Th}
and Proposition~\ref{rsksim} in~\cite{Kn}.
Proposition~\ref{rskeva} appears in \cite[Thm. D]{Kn2},
\cite[App. A]{St} for standard tableaux, and in
\cite[App. A1]{Fu} for tableaux of arbitrary weight.

Tableau switching, defined by means of the Bender-Knuth
transformations, was used by James and Kerber \cite[\S 2.8]{JK}
in a proof of the LR-rule.
More recently, the tableau switching map was defined in a more
general context in~\cite{BSS}; there its main properties were
established.
This is what we call the Tableau Switching map.
See also \cite[\S 2]{L2} for applications of tableau switching.
Proposition~\ref{switch} is Proposition~2.6 in \cite{BSS},
Proposition~\ref{swijeu} is contained in~\cite[\S 2]{BSS},
and Proposition~\ref{swilr} is a part of Theorem~3.1
in~\cite{BSS}.

The Littlewood-Robinson map (and other bijections of this type)
remains little studied despite being one of the oldest
Young tableau bijections.
The first such map was defined by Robinson~\cite{Ro} in
a different language in his effort to prove the LR-rule.
His proof was reworked by Macdonald in the first edition
of~\cite[I9]{Ma}.
For a detailed account on Robinson's paper see~\cite{L2};
in Section~2.5 of this paper, van Leeuwen defines another
bijection of Littlewood-Robinson type using the Tableau
Switching map.  This bijection, given here in Section~3.3,
is the one we call the Littlewood-Robinson map.
In fact, van Leeuwen proves in~\cite{L2} that
this map coincides with the original one defined
by Robinson and reworked by Macdonald.
From the point of view of the LR-rule this is unimportant,
since any such bijection yields a proof of the LR-rule.
Moreover, van Leeuwen shows using his definition that the
Littlewood-Robinson map enjoys nice properties, such as
Proposition~\ref{lirojeu} in this paper and
Corollary~2.5.2 in~\cite{L2}.
Proposition~\ref{lirojeu} follows from Proposition~\ref{swijeu},
and Proposition~\ref{liro} follows from Proposition~\ref{switch}
and the relations~(\maltese).

The problem of finding what we call a Fundamental Symmetry map
appears to be part of the `folklore' of the area;
it is very natural and has been considered independently
by several investigators (see e.g.~\cite{A1,A2,BZ,HK2}).
Proposition~\ref{funsim1} regarding our first map~$\rho_1$
is contained in~\cite[Example~3.6]{BSS}.
The Fundamental Symmetry maps $\rho_2$ and ${\rho_2}^\pr$ are new
to our best knowledge.
They are a byproduct of~\cite{PV}, and
were motivated by Fulton's appendix to~\cite{Bu}.
More precisely, $\ga(A)^\bullet$ is the composition of the linear
map $\Phi_l$, between LR-triangles and hives, defined
in~\cite[\S 4]{PV}, with Fulton's map in~\cite{Bu},
or equivalently the composition of the linear map
$\Psi_{l+1}\circ\Phi_{l+1}$,
between LR-triangles with last row equal to zero and
Berenstein-Zelevinsky triangles, defined in~\cite[\S 5]{PV},
with Carr\'e's map in~\cite[\S 3]{Ca}.
Thus, Proposition~\ref{carful} follows from~\cite{PV} and
Fulton~\cite{Fu} (which itself is based on Carr\'e's work~\cite{Ca}).
Note that both Carr\'e and Fulton's papers use
set of tableaux~$\FV^\ast(\la,\mu,\nu)$ to give combinatorial
interpretations of LR-coefficients~$c^\la_{\mu,\nu}$, which
they connect to BZ-triangles and hives, respectively.
In fact, the linear map~$\ga$ gives a simple combinatorial
proof of $c^\la_{\mu,\nu} = |\FV^\ast(\la,\mu,\nu)|$; in this
form it is new to the best of our knowledge (cf.~\cite{PV}).

The tableau $\tau(A)$ is called the companion tableau of $A$
in~\cite[\S 1.4]{L2}.
Proposition~\ref{tau} is equivalent to Proposition~1.4.3 there;
its proof is straightforward.
Proposition~\ref{funsim2} follows from Propositions~\ref{evajeu},
\ref{carful} and~\ref{tau}.

It is perhaps interesting to observe that the map $\Psi_l \circ \Phi_l$
from LR-triangles to BZ-triangles, given in~\cite[\S5]{PV}, essentially
contains $\gamma$ and $\tau$.
More precisely, let $A$ be a LR-tableau with LR-triangle $(a_{i,j})$,
and let $\Psi_l\circ \Phi_l (a_{i,j}) =(x_{i,j}, y_{i,j}, z_{i,j})$.
Note that the recording matrix $C=(c_{i,j})$ of $A$ satisfies
$c_{i,j} = a_{j,i}$.
Suppose $\gamma(A)= (d_{i,j})$ and $\tau(A)=(e_{i,j})$,
then $d_{i,j}=x_{l-j+1,\, l-j+i}$ and $e_{i,j}= y_{i,j-1}$ for all
$i<j$.
Besides, the numbers $d_{i,i}$ can be recovered from the $d_{i,j}$'s
and the $\mu_j$'s, and the numbers $e_{i,i}$ can be recovered
from the $e_{i,j}$'s and the $\nu_j$'s.

Finally, let us note that while $\rho_1$ is an {\em involution},
we do not know whether the same holds for~$\rho_2$ and ${\rho_2}^\pr$,
or equivalently, whether $\rho_2 ={\rho_2}^\pr$.
We suggest this in Conjecture~1 in Section~\ref{sec:further-alz}.

\bigskip
\subsection{Proof of relation ($\circledast$)}
\label{sec:proof-rel}


Since, by relation $(\maltese)$ one has that
$t_{l,\, k} = t_{k,\, l}^{-1}$, it is enough to prove that
\[
z_k \, t_{k,\, l}^{-1} \,z_l \, = \, z_{k+l} \, .
\]
This identity follows easily from the relations~$(\lozenge)$ of
the BK-transformations:
\[
\aligned
& z_k \, t_{k,\, l}^{-1} \,z_l \, = \,
[ (s_1)(s_2s_1)\cdots
(s_{k-1}s_{k-2}\cdots s_2s_1) ] \cdot
[ (s_{k}s_{k-1}\cdots s_2s_1)(s_{k+1}s_{k} \cdots s_3s_2)\\
&  \ \ \ \
(s_{k+2}s_{k+1}\cdots s_4s_3)\cdots
(s_{k+\ell-1}s_{k+\ell-2}\cdots s_{\ell+1}s_\ell) ] \cdot
[ (s_1)(s_2s_1)\cdots
(s_{\ell-1}s_{\ell-2}\cdots s_2s_1) ] \, \\
& = \,
[(s_1)(s_2s_1)\ldots
(s_{k-1}s_{k-2}\cdots s_2s_1)] \cdot
[(s_{k}s_{k-1}\cdots s_2s_1)\, (s_{k+1}s_{k}\cdots s_3s_2)(s_1)\\
&  \ \ \ \
(s_{k+2}s_{k+1}\cdots s_4s_3)(s_2s_1)\cdots
(s_{k+\ell-1}s_{k+\ell-2}\cdots s_{\ell+1}s_\ell)
(s_{\ell-1}s_{\ell-2}\cdots s_2s_1) ]\\
& = \, (s_1)(s_2s_1)\cdots (s_{k}s_{k-1}\cdots s_2s_1) \cdots
(s_{k+\ell-1}s_{k+\ell-2}\cdots s_2s_1) \, = \, z_{k+l} \,.
\endaligned
\]
Here the second equality follows from commuting
parenthesized products in $z_l$ with the the previous
products in~$t_{k,\, l}^{-1}$.

\bigskip
\subsection{Proof of lemmas}
\label{sec:proofs-lemmas}

\subsubsection{Proof of Lemma~1}
\label{sec:proofs-lemma1}

Let $k=\ell(B)$ and $l=\ell(A)$.
Use Propositions~\ref{evac} and~\ref{switch} to write Sch\"utzenberger
involution and Tableau Switching up to relabelling as products of
BK-transformations.
In the notation of Section~\ref{sec:prop-bk},
we need to prove that
$$t_{k,\, l} \,=  \, z_l \, z_{k+l}\, z_k .$$
This follows from relations $(\maltese)$ and $(\circledast)$.

\subsubsection{Proof of Lemma~2}
\label{sec:proofs-lemma2}

Let $A\in\YT(\la, \ba)$ with recording matrix $V=(v_{i,j})$.
Then, by Proposition~\ref{rsk} for $V^\updownarrow=(v_{k+1-i,\, j})$,
we have $\vp(V^\updownarrow)=(A, -)$.
Since ${V^\updownarrow}^\ast =(v_{i,\, k+1-j}) =U$,
Proposition~\ref{rskeva} implies $\vp (U)= (\xi^\N(A), -)$, and
the claim follows.

\subsubsection{Proof of Lemma~3}
\label{sec:proofs-lemma3}

First, we need to show that~$\wh A$ and $\wh B$ are
both LR-tableaux.  Indeed, since~$C$ is canonical, the integer
$(i)$~appears $\alpha_i$ times in~$C$, which is the total number
of times~$(i+1)$ appears in~$A$. Therefore, when reading the
$\word(\wh A)$, the number of integers~$(i)$ is always at least
as many as the number of integers~$(i+1)$.
By definition, this implies that $\wh A$ is a LR-tableau,
and the same argument works for~$\wh B$.

Now observe that tableaux~$C$ and~$D$ remain unchanged
under~$\zeta$.  Furthermore, since $\zeta^\LRN$ is applicable
and well defined, we easily see that action of~$\zeta^\LRN$
restricted to $(B,A)$ coincides with the action of~$\zeta$
on~$(B,A)$.  Indeed, simply observe that BK-transformations
commute with taking compositions of tableaux.  Thus, so do
elements~$t_{r,\, m-r}$ and by Proposition~2, the Tableau Switching
map~$\zeta$.  Therefore, the restriction of~$\zeta^\LRN$
to~$(B,A)$  coincides with~$\zeta$, which implies the result.

\subsubsection{Proof of Lemma~4}
\label{sec:proofs-lemma4}


We need the following property of tableau switching
given in~\cite[Thm.~2.3]{BSS}.
Let $(X,Y) = \zeta(U, V \star W)$; then, there is an alternative
way to calculate $(X,Y)$.
Start with $(V^\pr, U^\pr) = \zeta(U,V)$, and compute
$(W^\pr,U^{\pr\pr}) = \zeta(U^\pr,W)$.
Then $(X,Y) = (V^\pr\star W^\pr, U^{\pr\pr})$.
By the symmetry (\maltese), the same ``distributivity'' property
holds for $(X^\pr,Y^\pr) = \zeta(U \star V, W)$.
During the proof we will adopt the following convention:
if~$U$ and~$V$ are tableaux filled with integers $1,\dots, k$
then by $U\star V$ we will denote the tableau obtained by
relabelling the entries of $V$ by adding $k$ to them, and then
taking the union of $U$ and $V$ (if this is possible).

Let $(\overline{A},\overline{B})=\zeta(B,A)$, and let $A^\pr$, $B^\pr$
be defined as in reduction~\ref{sec:zoo-462}.
We have to show that $\overline{A}=A^\pr$ and $\overline{B}=B^\pr$.
Now, by Proposition~\ref{funsim1}, the map~$\rho_1$ is a special case
of the Tableau Switching map~$\zeta$.
The first application of~$\zeta$ computes
$(\Can(\tau),C) = \zeta(\Can(\mu),B)$.
The second application of~$\zeta$ computes
$(\Can(\vk),\wh B^\pr)$ as the switching of
$(\Can(\wh \tau),(C\comp_{\la_1,0} G)\star A)$.
By decomposing $\Can(\wh \tau)$ as $ \Can(\de) \star \Can( \tau)$,
and $\Can(\vk)$ as $\Can(\widetilde{\mu}) \star \Can(\nu)$ we
obtain
\begin{equation}
(\Can(\widetilde{\mu}) \star\Can(\nu), F\star B^\pr)=
\zeta(\Can(\de)\star\Can(\tau), (C\circ_{\la_1,0}G)\star A)\, .
\tag{$\blacklozenge$}
\end{equation}

Recall that tableau switching commutes with taking compositions
(see Subsection~\ref{sec:proofs-lemma3} above), and observe that
$\Can(\tau)$ lies to the left and below~$G$.
Now, by the ``distributivity'' property,
the second application of~$\zeta$
starts with tableau switching of~$\Can(\tau)$ and~$C$,
which is the inverse of the first application of~$\zeta$.
The resulting tableau~$B$ is then switched with~$A$, giving~$B^\pr$
as desired.  The remaining steps to be done according to the
``distributivity'' property as above do not change the
tableaux~$B^\pr$ as it contains the largest integers~$k+1,\ldots,2k$.
Therefore, $\overline{B}=B^\pr$.
The third application of $\zeta$ switches~$F$ with
$\Can(\widetilde{\mu})\star\Can(\nu)$, but according to
$(\blacklozenge)$, $\Can(\tau)$ switches first with
$(C\circ_{\la_1,0} G)\star A$ yielding
$(\Can(\mu)\circ_{\la_1,0} G)\star \overline{A}, \overline{B})$;
then $\Can(\de)$ switches with
$(\Can(\mu)\circ_{\la_1,0} G)\star \overline{A}$ yielding
$(\Can(\widetilde{\mu})\star\Can(\nu), F)$.  Therefore
$\rho_1(F) = (\Can(\mu)\circ_{\la_1,0} G)\star \overline{A}$,
and restricting this tableau to the last $k$ integers gives
$\overline{A} =A^\pr$, as desired.

It remains to show that $\rho_1$ is applicable the
three times in the circuit.
Since~$B$ is already LR-tableau, the first application is valid.
By Proposition~\ref{funsim1}, the resulting tableau~$C$ is also~LR.
We need to show that $(C \comp_{\la_1,0} G)\star A$
is LR-tableau.
Since~$(C \comp_{\la_1,0} G)$ and~$A$ are already LR-tableau,
all we need to show is that the number of $(k+1)$'s is always
at most the number of $k$'s in a word.  But that is clear since
there are $s=\nu_1$ integers~$k$ in~$G$, which all appear before
the word reaches~$A$.
Finally, since $F$ is obtained by switching from
$\Can(\de)$, it is a LR-tableau by Proposition~\ref{swilr}.
This justifies the third application of~$\rho_1$ and
completes the proof of the lemma.

\subsubsection{Proof of Lemma~5}
\label{sec:proofs-lemma5}


Let $A\in \YT(\la/\mu, \ba)$.
We have to show that $A^{\bw} = \chi(A)$.
For this we will use the following way of computing $\xi^\N(D)$
for a tableau $D$:
Write $D=E\star F$, where $E$ has integers $1,\dots, l$, and
$F$ has integers $l+1,\dots, l+k$.
First compute $\xi(E)$; then relabel the entries of $F$
by subtracting $l$ from them,
and compute $(F^\pr, E^\circ) =\zeta(\xi(E), F) $.
Finally, compute $F^\circ= \xi(F^\pr)$, and relabel the entries
of $E^\circ$ by adding $k$ to them.
We obtain~$\xi(D) = F^\circ \star E^\circ$.
This follows immediately from relations~$(\circledast)$ in
Section~\ref{sec:prop-bk}.
Now $B^\circ =\xi^\N(B)$ is computed as follows:
Observe that $ \xi^\N(C^\circ) = C$, let $(A^\pr, C^\pr)=\zeta(C,A)$.
Then, up to relabelling of $C^\pr$, we have
$B^\circ=\xi^\N(A^\pr)\star C^\pr$.
By relation~$(\maltese)$, we have $A^{\circ\circ}=A^\pr$.
Thus~$B^{\circ\circ}=A^\pr\star C^\pr$.
Finally, $B^{\bw}$ is obtained by taking
$(C^{\pr\pr}, A^{\pr\pr})=\zeta(\xi^\N(A^\pr), C^\pr)$,
and computing $\xi^\N(C^{\pr\pr})$.
Since $C^{\pr\pr} = C$, $\xi^\N(C^{\pr\pr}) = C^\circ$.
Thus, up to relabelling, $A^{\bw}=A^{\pr\pr}$.
We claim that $A^{\pr\pr}=\chi(A)$.
Recall that by definition of the reversal map~$\chi$,
the image~$\chi(A)$ is the second component of
$\zeta( \xi(A^\pr), C^\pr)$.  This implies the result.

\subsubsection{Proof of Lemma~6}
\label{sec:proofs-lemma6}
Let $\mu = (\mu_1,\ldots,\mu_\ell)$, $k:= \ell+1$, and
set~$\mu_k = 0$.  Compute $a_{1,r} := \mu_{k-r} - \mu_{k-r+1}$,
for all $1 \le r \le \ell$, and
$a_{i,j} := a_{1,j-i+1}$, for all $1 \le i \le j \le \ell$.
It is easy to see that $(a_{i,j})$ is the recording matrix of
the desired tableau $A = \xi(\Can(\mu))$.

\subsection{Proof of theorems}
\label{sec:proofs-thms}

\subsubsection{Proof of Theorem 1}
\label{sec:proofs-thm1}

\begin{figure}[hbt]
\vskip.3cm
\psfrag{a}{$\vp$}
\psfrag{b}{$\psi$}
\psfrag{c}{$\phi$}
\psfrag{d}{$\zeta^\N$}
\psfrag{e}{$\xi^\N$}
\psfrag{f}{$\rho_1$}
\psfrag{g}{$\zeta^\LRN$}
\psfrag{h}{$\zeta$}
\psfrag{x}{$\rho_2$}
\psfrag{y}{$\chi$}
\psfrag{1}{$1$}
\psfrag{2}{$2$}
\psfrag{3}{$3$}
\begin{center}
\epsfig{file=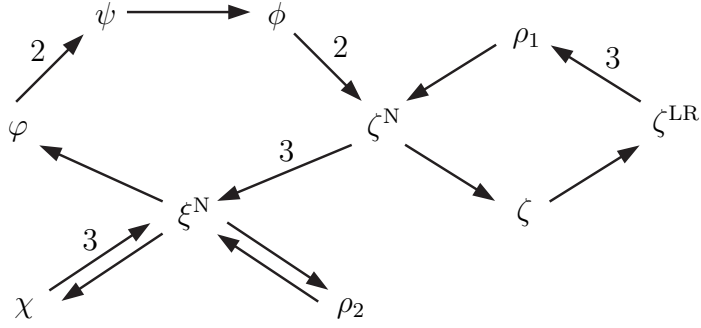,width=9.2cm}
\end{center}
\caption{Diagram of linear reductions.
}
\label{f:graph}
\end{figure}

As we showed in Section~\ref{sec:zoo}
(see Section~\ref{sec:zoo-together} and Corollary~3),
all the maps listed
in Theorem~1 are linearly equivalent.  To prove the second
part of the theorem, let us summarize all linear reductions
in Figure~\ref{f:graph}.  Here we draw an arrow for every
linear reduction given by a ps-circuit~$\CC$, and place the
cost of the circuit~$\bs(\CC)$ above the arrow.
We do not to write the cost above trivial (cost~1) circuits.


Recall the second part of Composition Lemma which claims that
one needs to take a product of costs when taking a composition
of linear reductions.  Now observe that from each map in the
diagram one can go to any other map (taking arrows in the
reverse direction) so that the product never exceeds~$36$.
This product maximizes when going from~$\rho_1$ to~$\chi$.
The verification is straightforward and left to the reader.

\subsubsection{Proof of Theorem 2}
\label{sec:proofs-thm2}

Consider the map~$\xi$.  By Proposition~1, it can be computed
at a cost of $\binom{k}{2}$ BK-transformations.  Each
BK-transformation is a piecewise linear map which can be
computed at a cost of~$O(k)$ additions and $\max/\min$
operations on integers~$a_{i,j}$.
Note that the size of integers~$a_{i,j}$ never exceed~$\la_1$,
so during and after all BK-transformations they have
bit-size~$O(\log m)$.  Therefore,  the total cost of
computing~$\xi$ is~$O(k^3 \log m)$, as in the theorem.
Also, by definition~$\xi$ is a size-neutral map
(see Section~\ref{sec:basic-comp}).


By Theorem~1, all other maps are linearly reducible to~$\xi^\N$.
Denote by~$\al$ any of the remaining maps in Theorem~2, and
by~$\CC$ denote a $\xi^\N$-based ps-circuit computing~$\al$
(at a cost at most~$36$). Recall that~$\xi^\N$ and all linear
cost maps are size-neutral (see Section~\ref{sec:main-red}),
which makes map~$\al$ size-neutral as well.
Thus, the cost of computing
any of the (at most~$36$) maps~$\xi^\N$ is~$O(k^3 \log m)$,
where parameters~$k$ and~$\log m$ are linear in the input
parameters as in Theorem~2.  Therefore, the total cost of
computing~$\al$ following~$\CC$ is~$O(k^3 \log m)$, as
desired.

\bigskip


\section{Further bijections} \label{sec:further}

\subsection{Third Fundamental Symmetry map}
\label{sec:further-alz}

Here we present another Fundamental Symmetry map~$\rho_3$,
defined in~\cite{A1,A2}.  The definition in these papers
is somewhat convoluted so we restate it here for the sake
of clarity.

\smallskip

Start with LR-tableau $A \in \LR(\la/\mu,\nu)$.
Fill shape~$[\mu]$ with zeros.  We remove rows one by one,
beginning with the bottom row.  In each row to be removed,
build a chain of integers in previous rows, starting with
the last element and going to the first element.  For each
such element~$x$, find the largest element~$y <x$ in the
previous row, not used by the previous chains (starting
from row containing~$x$), then the largest element $z <y$
in the row above that of~$y$ not used by the previous
chains, etc.  Now replace~$y$ with~$x$,~$z$ with~$y$, etc.
unless the integer~$k$ goes in $<k$-th row; stay put in
that case.  Note that each zero forms a chain of length~1.

Denote by $v_{i,j}$ the number of chains of length~$(i-j+1)$
which start from~$i$-th row.  Let~$B \in \YT(\la/\nu,\mu)$
be a Young tableau corresponding to recording
matrix~$V = (v_{i,j})$. We claim that~$B$ is a LR-tableau,
and define $B = \rho_3(A)$.

\bigskip

\begin{figure}[hbt]
\vskip.3cm
\psfrag{A}{$A$}
\psfrag{B}{$B$}
\psfrag{V}{$V = $}
\psfrag{r}{$\rho_3$}
\begin{center}
\epsfig{file=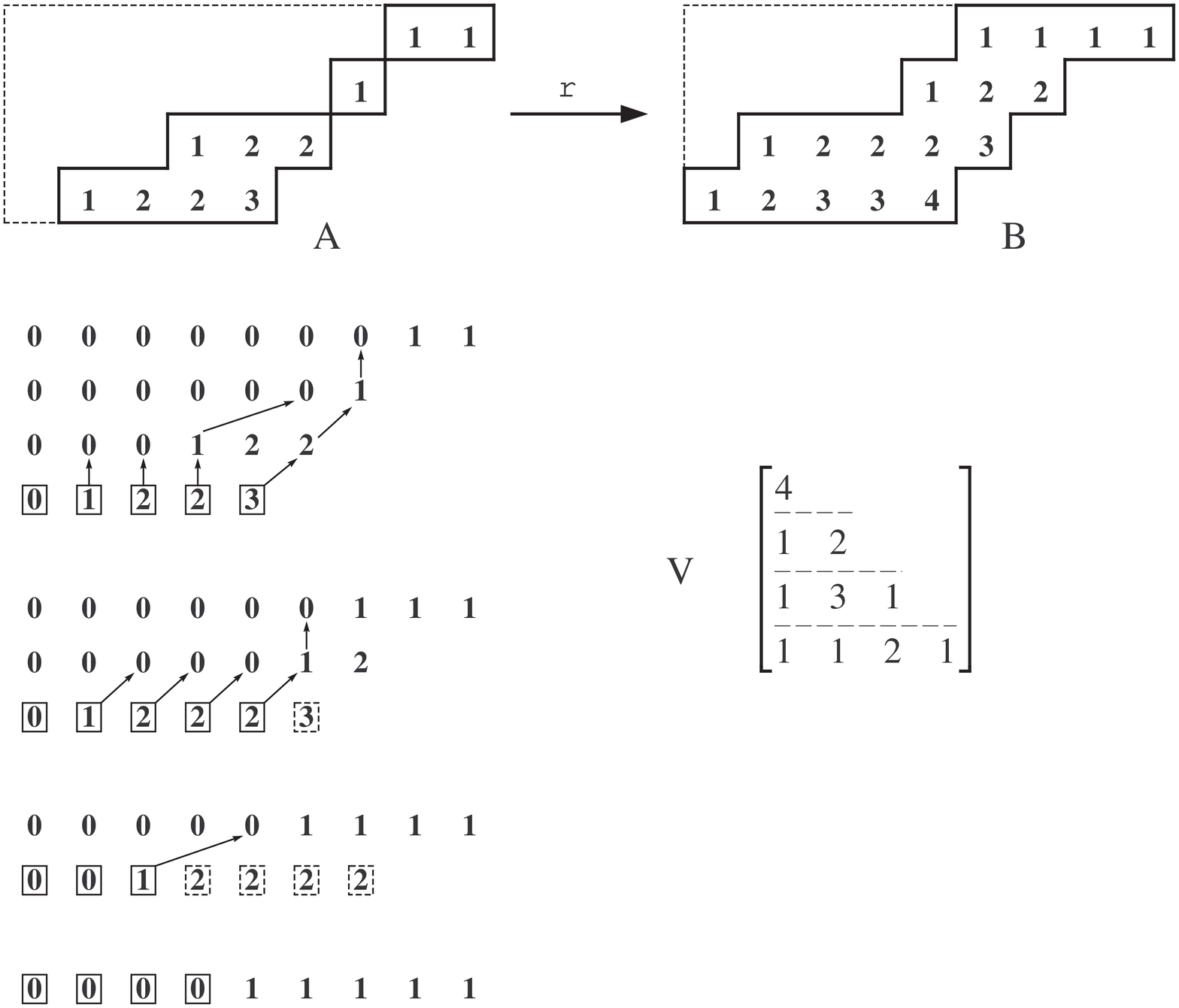,width=13.3cm}
\end{center}
\caption{An example of the map $\rho_3: A \to B$, where
$A \in \LR(\la/\mu,\nu)$, $B \in \LR(\la/\nu,\mu)$, and
$\la = (9,7,6,5)$, $\mu = (7,6,3,1)$, $\nu = (5,4,1)$.
There is one chain of length~$4$, $3$, $1$, and two chains
of length~$2$, starting from the $4$-th row
(see $4$-th row of~$V$).
}
\label{f:alz}
\end{figure}

\bigskip

{\bf Conjecture~1.} \ {\it The Fundamental Symmetry maps
$\rho_1$, $\rho_2$, ${\rho_2}^\pr$ and $\rho_3$  are identical.}

\bigskip

The conjecture is supported by numerical evidence.
Also, in~\cite[\S 5]{A2} it was shown that $(\rho_3)^2=1$ by an
involved argument.
If established, the conjecture would simplify the proof
of Theorem~1  and further emphasize the importance of
the fundamental symmetry.

\bigskip

\subsection{Inverse maps}
\label{sec:further-inv}

Recall that maps~$\vp$ and~$\phi$ are one-to-one.

\bigskip

{\bf Conjecture~2.} \ {\it The RSK map~$\vp$ and
the Littlewood-Robinson map~$\phi$ are linearly
equivalent to their inverses. }

\bigskip

Let us emphasize here the need to distinguish
between direct and inverse maps.  It is a well
known and studied phenomenon in Cryptography
that some maps are easily computed, while their
inverses are not; taking powers over the finite field
vs. taking discrete logarithm being the most celebrated
example.  The conjecture above says there is no such problem
with Young tableau maps and  all Young tableau bijections
are linearly equivalent to their inverses.  Note that
Sch\"utzenberger involution, Tableau Switching, Reversal
and the First Fundamental Symmetry maps are {\em equal}
to their inverses~$(\maltese)$, and if Conjecture~1 holds
so is the Second Fundamental Symmetry map.  Thus the problem
makes sense only for RSK and Littlewood-Robinson maps, as in
the conjecture.

\bigskip

\subsection{Octahedral map}
\label{sec:further-oct}

Suppose four partitions $\la, \mu, \nu$, and $\tau$ satisfy
$|\tau| = |\la| + |\mu| + |\nu|$.  The \emph{Octahedral map}
is a one-to-one correspondence:
$$\varsigma: \, \bigcup_{\si \, \vdash |\la|+|\mu|} \,
\LR(\si/\la,\mu)
\times \LR(\tau/\si,\nu) \ \too
\bigcup_{\pi \, \vdash |\mu|+|\nu|} \,
\LR(\pi/\mu,\nu)
\times \LR(\tau/\la,\pi).
$$
A bijection of this type was introduced in~\cite{KTW}, in the
equivalent language of hives, as a tool for a new proof of the
LR-rule. This is defined using an {\em octahedral recurrence}
considered earlier in connection with the enumeration of
alternating sign matrices~\cite{RR}, and has recently
appeared in other context~\cite{HK2,Sp,Thu}.
%


As it was the case with the Littlewood-Robinson map, the existence
of any such bijection~$\vs$ will suffice in the proof of the
LR-rule given in~\cite{KTW}.
So, we will define an alternative version of the octahedral map
using the tableau switching map~$\zeta$.
From our point of view, this has the advantage that the definition
given here yields the reduction $\vs \lrm \zeta$.

We define an Octahedral map $\vs$ as follows:
Let $A\in\LR(\si/\la,\mu)$, $B\in\LR(\tau/\si,\nu)$.
Consider $(A^\pr, C^\pr)=\zeta(\Can(\la), A)$; since $A$
is a LR tableau, $A^\pr=\Can(\mu)$ and
$C^\pr \in \LR(\si/\mu,\la)$.
Next, let $(B^\pr, C^{\pr\pr})=\zeta(C^\pr, B)$.
Again, since switching preserves the property of being a
Littlewood-Richardson tableau, there is some partition
$\pi$ of $|\mu| + |\nu|$, such that $B\in \LR(\pi/\mu, \nu)$
and $C^{\pr\pr} \in \LR(\tau/\pi, \la)$.
Finally, let $(C^{\pr\pr\pr}, D)=\zeta(\Can(\pi), C^{\pr\pr})$.
Thus, $C^{\pr\pr\pr} =\Can(\la)$ and $D\in\LR(\tau/\la,\pi)$.
We define $\vs(A,B)=(B^\pr, D)$.

The following proposition is a consequence of $(\maltese)$.

\begin{prop}
The map $\vs$ defined above is a bijection.
\label{octa}
\end{prop}

\medskip
\noindent
{\bf Corollary~4.} \ {\it The Octahedral map~$\vs$ is
linearly reducible to Tableau Switching map~$\zeta$ and
other maps in Theorem~1.}

\medskip

Note that~$\vs$ is defined by a simple sequential circuit
which uses the map $\zeta$ three times; thus $\vs \lrm \zeta$.
Another way to prove this is to show first that~$\vs$ is
a composition of the $\LR$-Tableau Switching map~$\zeta^\LRN$
and the fundamental symmetry maps.
Let us conclude this section with the following natural
conjecture:

\bigskip
\noindent
{\bf Conjecture~3.} \ {\it The map~$\vs$ and the map defined
in~\cite{KTW} are identical and linearly equivalent to maps
in Theorem~1. }


\medskip

We should mention that the connection between Tableau
Switching map~$\zeta$ and the Octahedral map in~\cite{KTW}
follows from~\cite{HK1} through linear equivalence with the
Jeu de Taquin map~$\psi$. Also, it was shown in~\cite{HK2}
that in a special case the map~$\vs$ gives a (another version)
fundamental symmetry map (see the ``commutor'' in Section~5.2
of~\cite{HK2}).  It is natural to conjecture that this
fundamental symmetry map coincides with~$\rho_1$ as well.



\bigskip
\subsection{Burge correspondence}
\label{sec:burge}


Let $\widetilde\vp$ denote the \emph{Burge
correspondence}~\cite{Burge} (see also~\cite[A4.1]{Fu}).
Numerically, it defines a one-to-one map between the same
sets as the RSK map~$\vp$~:
$$
\widetilde\vp: \, \Mat(\ba,\bb) \, \too \, \bigcup_{\la \vdash N} \,
\YT(\la,\bb) \times \YT(\la,\ba)
$$
This bijection is related to RSK correspondence in the following way.
Let $V=(v_{i,j})$; denote $V^\updownarrow :=(v_{k+1-i,j})$ and
$V^\leftrightarrow := (v_{i,k+1-j})$.
Let $\vp(V)=(B,A)$.
Then, since column insertion commutes with row insertion~\cite[A2]{Fu},
we have that $\widetilde\vp(V^\updownarrow) =(B, \ast)$ and
$\widetilde\vp(V^\leftrightarrow)=(\ast, A)$.
Thus $\vp \lrm \widetilde\vp$ and this is done by a $\widetilde\vp$-based
simple parallel circuit which uses $\widetilde\vp$ twice.
Similarly, one can show that $\widetilde\vp \lrm \vp$, which
implies~$\vp \sim \widetilde\vp$.

\subsection{Hillman-Grassl map}
\label{sec:further-hg}

Let $\la \vdash n$ be a fixed partition, $\ell = \ell(\la)$,
$m = \la_1$.  For every function $F: [\la] \to \zz_{\ge 0}$ and every
$-\ell< c < m$, define {\em diagonal sums}
$$\al_c(F) \, = \, \sum_{(i,j) \in [\la], \, j-i = c} \, F(i,j),
$$
and {\em rectangular sums}
$$\be_c(F) \, = \  \sum_{i=1}^{i_c} \ \sum_{j=1}^{j_c} \ F(i,j),
$$
where $(i_c,j_c)$ is the last square on the diagonal $j-i=c$.
Now, let $\bdd = (d_{1-\ell},\ldots,d_{m-1})$ be a nonnegative
integer array.  Define~$\cb_\bdd$ to be sets of all nonnegative integer
functions~$F$ as above, such that $\be_c(F) = d_c$,
for all~$-\ell< c < m$.  Similarly, define~$\ca_\bdd$
to be sets of all reverse plane partitions of shape~$\la$,
such that $\al_c(F) = d_c$, for all $-\ell< c < m$.


The {\em Hillman-Grassl} (HG) {\em bijection} defines a
one-to-one map
$\vt_\la: \cb_\bdd \to \ca_\bdd$ \cite{HG} (see also~\cite{Ga,GH,P2}).
It is easy to check that when~$\la = (k^k)$ the set~$\cb_\bdd$
coincides with $\Mat(\ba,\bb)$
for certain~$\ba,\bb$,  while~$\ca_\bdd$ corresponds to pairs
of GT-patterns joined at the diagonal~\cite{P2}.
In~\cite[Thm. 10.2]{Ga} Gansner showed  that the map
 $\vt_k:=\vt_\la$ in this case coincides, up to a linear cost map,
 with the Burge correspondence~$\wt\vp$.
 This immediately gives~$\vt_k \sim \widetilde\vp$.
Combining this equivalence with the one in the previous section
we conclude that the HG-map in the ``square case''
is linearly equivalent
to the RSK correspondence~$\vt_k \sim \vp$, as well as all other maps
listed in Theorem~1.  We omit the (easy) details.

For general shapes~$\la$, given~$F \in \ca_\bdd$,
set $k:=\max\{m,\ell\}$, and fill with zeros the rest
of the $k \times k$ square containing~$[\la]$.
Now apply RSK map~$\vp$ to the resulting matrix.
At the end, join at the diagonal the two GT-patterns of
the resulting tableaux, and restrict this function to squares
in~$[\la]$ (see~~\cite{P2} for details).
We leave it to the reader to show that this
defines linear reduction~$\vt_\la \lrm \vp$, proving linear
equivalence~$\vt_\la \sim \vp$ in general case, where~$\la$
is a part of the input.

To conclude, we note that the connection of~$\vt$ with
BK-transformations was observed in~\cite{P2},
where it was also shown that~$\vt$ can be computed
in~$O(k^3 \log m)$, where $m = \max_c\{d_c\}$.



\bigskip

\subsection{Other symmetry maps}
\label{sec:further-sym}

Beside fundamental symmetry maps, there is a large number of
``hidden'' symmetries of Littlewood-Richardson
coefficients~$c_{\mu,\nu}^\la$.  These symmetries form a finite
group and were studied on a number of occasions (see~\cite{BZ}).
As we mentioned
above, a subgroup of index~2 of these symmetries can be given
by linear cost maps~\cite{PV}.  Since the fundamental symmetry
is a remaining generator, the symmetries outside this subgroup
are given by maps which are all linearly equivalent to~$\rho_1$.

\smallskip

A different kind of symmetry map was given in~\cite{HS}
(see also~\cite[\S 3]{A0}):
$$\varrho: \, \LR(\la/\mu,\nu) \, \too \, \LR(\la^\pr/\mu^\pr,\nu^\pr),
$$
where $\la^\pr$ denotes the conjugate diagram
(reflected across $i=j$ line).  The bijections given
in~\cite{HS,A0} use a modified insertion map.
It would be interesting to see
whether this symmetry map is linearly equivalent to maps we
consider in Theorem~1.

\bigskip

\subsection{Sch\"utzenberger involution}
\label{sec:further-sch}

Recall that we have not been able to show that the general
Sch\"utzenberger involution~$\xi$ reduces
to the other bijections appearing in this paper, while we
were able to show that~$\xi^\N$ and~$\chi$ do.
If~$\xi$ were not reducible to the other bijections this would
mean that $\chi$ is a more natural extension of $\xi^\N$
to skew shapes than~$\xi$.
Recall \cite[\S 5]{BSS} that reversal is also more natural than
Sch\"utzenberger involution from the point of view of dual equivalence.
Proving that~$\xi$ is reducible to~$\xi^\N$ remains an open problem.
\bigskip

\section{Final Remarks}
\label{sec:final}

{\bf 1.} \ Note that we never attempted to give a
lower bound on the complexity of the cost of tableau
bijections.  Since all constructions require $\theta(k^2)$
min-max operations, it is conceivable that such lower
bound can be obtained
by means of {\it Algebraic Complexity Theory}~\cite{BCS}.
In other words, if one properly restrict
the class of algorithms to consider, the lower bound $\Omega(k^3)$
might be attainable.  Further investigation of this matter
would be of great interest.

\bigskip

\noindent
{\bf 2.}  \  While the Octahedral map defined in~\cite{KTW}
(see also~\cite{HK2}) looks extremely natural,
it lacks formal and complete treatment
in combinatorics literature.  Our alternative version
of this map and Conjecture~3 is a further indication of this.
%
%
We would like to encourage the reader to further study this
map and its connections to other combinatorial maps.

\bigskip

\noindent
{\bf 3.}  \ There are a score of other notable Young tableau
bijections not mentioned in the paper.  While some of these
are based on some kind of insertion/evacuation procedures and
thus seem strongly related to the maps we study, others are
of a different nature.
Examples of the first type include
Lascoux-Sch\"utzenberger action of~$S_m$ on~$\YT(\la/\mu; m)$,
RSK for shifted and super tableaux, etc.  Examples of the
second type include the Novelli-Pak-Stoyanovskii's bijection, and
nonintersecting paths arguments.  We refer to~\cite{Fu,Sa,St}
for definitions and references.  It would be nice to place these
maps into our framework and perhaps even introduce some kind
of complexity style hierarchy on them.

\bigskip

\noindent
{\bf 4.}  \ In~\cite{BSS} the authors present an important
characterization of the tableau switching through its natural
properties.  In principle, one can use our linear reductions
to obtain similar characterizations of the fundamental
symmetry maps.  We challenge the reader to make such a
characterization explicit.
This could lead to a positive resolution of Conjecture~1
and give a better understanding of the subject.

\bigskip

\noindent
{\bf 5.}  \  It was noted in~\cite{P2} that the Hillman-Grassl
map extends to real-valued functions and is given by a
continuous  piecewise-linear volume-preserving map in this case.
Similar observations were made for other sets of Young tableaux
and other tableau bijections (see e.g.~\cite{A1,BZ,KB,PV}).
One can check that most of our linear reductions also
extend to real-valued tableaux, and our linear cost maps are
in fact (the usual, geometric) linear maps.
Thus, one ``real extension''
suffice to establish the others.  We leave further exploration
of this subject to the interested reader.

\bigskip

\noindent
{\bf 6.}  \ We make no effort to optimize the constant~36
in the proof of Theorem~1.  In fact, if $\rho_1=\rho_2$
as Conjecture~1 suggests, this constant immediately
drops down to~12.  We wanted to emphasize the existence
of such a constant, and would be just as happy if it worked
out to be around~200, as in the first draft of the paper.

The idea to study the cost as the number of ``large operations''
is well-known in Computer Science literature.  It was rejected
for the purposes of Cryptography after A. Shamir showed in~\cite{Sh}
that factoring has a polynomial cost algorithm in the number
of arithmetic operations (of possibly very large integers).
In our situation, all maps are size-neutral, so this problem
never arises.  Still, we should warn the reader that the
constant implied by the $O(\cdot)$ notation in Theorem~2 is
much larger than~36.

\bigskip

\noindent
{\bf 7.}  \ It was noticed in~\cite{P3} that there is a certain
uniqueness behind the linear cost partition bijections.  This
method was later used in~\cite{PV} to obtain several linear cost
LR-symmetry maps.  We wonder if there is an ``automatic'' procedure
to define the maps in Theorem~1.  Even a framework
for that would be of great interest.

\bigskip

\noindent
{\bf 8.}  \ We conclude by saying that even though we were
able to formalize the connections between Young tableau
bijections, this is important on a computational and perhaps
philosophical level, but is hardly an ``explanation from the
Book".  We firmly believe that the real underlying structure
behind these connections lies in the study of canonical bases
in representation theory of symmetric and full linear groups.

\bigskip



\vskip 1.3pc
{\bf Acknowledgements}

\smallskip
\noindent
We are grateful to Olga Azenhas, Arkady Berenstein,
Andre Henriques, Leonid Levin, Ezra Miller,
Richard Stanley and Terry Tao
for interesting conversations and helpful remarks.
Special thanks to David Jackson for reading and editing
the previous version of the paper and his meticulous
attention to details.


The first author was supported by the NSA and the NSF.
Part of this work was done during second author's sabbatical
stay at MIT Mathematics Department.
The second author would like to thank CONACYT and DGAPA-UNAM for
financial support, and to Richard Stanley for his support in
the realization of the visit.
The second author was also partially supported by DGAPA-UNAM
IN111203.  The first author would also like to thank the second
author for organizing his visit to Mexico, where this work
was continued.

\bigskip




{\small
\begin{thebibliography}{99}

\bibitem{A0} O.~Azenhas,
The admissible interval for the invariant factors
of a product of matrices,
{\em Linear and Multilinear Algebra} {\bf 46} (1999), 51--99.

\bibitem{A1} O.~Azenhas,
Littlewood-Richardson fillings and their symmetries,
in {\em Matrices and Group Representations},
Textos de Matem\'atica, S\'erie B {\bf 19} (1999),
Depto. de Matem\'atica, Univ. de Coimbra, 81--92.

\bibitem{A2} O.~Azenhas, On an involution on the set of
Littewood-Richardson tableaux and the hidden commutativity,
preprint, 2000,
available from
{\tt http://dingo.mat.uc.pt/}$\widetilde\, $
{\tt cmuc/publicline.php?lid=1}

\bibitem{BF}
A. Barvinok  and S. Fomin,
Sparse interpolation of symmetric polynomials,
{\em Adv. Appl. Math.} {\bf 18} (1997),  271--285.

\bibitem{BK}
E.~A. Bender and D.~E. Knuth,
Enumeration of plane partitions,
{\em J. Combin. Theory}, Ser. A  {\bf  13} (1972), 40--54.

\bibitem{BSS} G.~Benkart, F.~Sottile and J.~Stroomer,
Tableau switching: algorithms and applications,
{\em J. Combin. Theory}, Ser. A   {\bf 76} (1996), 11--43.


\bibitem{BZ} A.~D. Berenstein and A. Zelevinsky,
Triple multiplicities for $sl(r+1)$ and the spectrum of the exterior
algebra of the adjoint representation,
{\em  J. Algebraic Combin.} {\bf 1} (1992), 7--22.

\bibitem{Bu} A. Buch,
The saturation conjecture (after A. Knutson and T. Tao),
{\em Enseign. Math.} {\bf 46} (2000), 43--60.

\bibitem{Burge}
W.~H. Burge,
Four correspondences between graphs and generalized Young tableaux,
{\em J. Combin. Theory}, Ser. A  {\bf 17} (1974), 12--30.

\bibitem{BCS}
P. B\"urgisser, M. Clausen and M.~A. Shokrollahi,
{\em Algebraic Complexity Theory},
Grund\-lehren der mathematischen Wissenschaften,
Vol. 315, Springer, Berlin, 1997.


\bibitem{Ca} C. Carr\'e,
The rule of \liri in a construction of Berenstein-Zelevinsky,
{\em Internat. J. Algebra \& Comput.} {\bf 1} (1991), 473--491.

\bibitem{Fu} W. Fulton,
{\em Young Tableaux},
L.M.S. Student Texts 35,
Cambridge U. Press, 1997.

\bibitem{Ga}
E.~R. Gansner,
Matrix correspondences of plane partitions,
{\em Pacific J. Math.} {\bf 92} (1981), 295--315.

\bibitem{GJ}
M.~R. Garey  and D.~S. Johnson,
{\em Computers and intractability:
A guide to the theory of $\textsc{NP}$-completeness},
Freeman, New York, 1979.


\bibitem{GH}
R.~M. Grassl and A.~P. Hillman,
Functions on tableau frames,
{\em Discrete Math.} {\bf 25} (1979), 245--255.

\bibitem{Ha}
M.~D. Haiman,  Dual equivalence with applications,
including a conjecture of Proctor,
{\em Discrete Math.} \textbf{99} (1992),  79--113.

\bibitem{HS}
P. Hanlon and S. Sundaram,
On a bijection between Littlewood-Richardson fillings of
conjugate shape
{\em J. Combin. Theory}, Ser. A  {\bf 60} (1992), 1--18.

\bibitem{HK1} A. Henriques and J. Kamnitzer,
On the symmetric monoidal category constructed using the
octahedral recurrence, preprint, 2003.

\bibitem{HK2} A. Henriques and J. Kamnitzer,
The octahedron recurrence and $gL_n$ crystals,
{\em Proc. $16$-th FPSAC Conf.}, Vancouver, BC, 2004,
available from {\tt
http://www-math.mit.edu/} $\widetilde\,$~{\tt andrhenr/}.

\bibitem{HG}
A.~P. Hillman and R.~M. Grassl,
Reverse plane partitions and tableau hook numbers,
{\em J. Combin. Theory}, Ser. A  {\bf 21} (1976), 216--221.

\bibitem{JK}
G.~D. James and A. Kerber,
{\em The Representation Theory of the Symmetric Group},
Encyclopedia of Mathematics and its Applications {\bf 16},
Addison-Wesley, Reading, Massachusetts, 1981.

\bibitem{KB}
A.~N. Kirillov and A.~D. Berenstein,
Groups generated by involutions, Gelfand-Tsetlin patterns,
and combinatorics of Young tableaux,
{\em St. Petersburg Math. J.} {\bf 7} (1996), 77--127.

\bibitem{Kn} D.~E. Knuth,
Permutations, matrices and generalized
young tableaux, {\em Pacific J. Math.} {\bf 34} (1970),
709--727.

\bibitem{Kn2} D.~E. Knuth,
{\em The Art of Computer Programming}, vol.~3,
{\em Sorting and Searching},
second ed., Addison-Wesley, Reading, Massachusetts, 1998.

\bibitem{KTW} A. Knutson, T. Tao and C.~T. Woodward,
A positive proof of the Littlewood-Richardson rule using the
        octahedron recurrence, arXiv:math.CO/0306274.

\bibitem{L1}
 M.~A.~A. van  Leeuwen,
The Robinson-Schensted and Sch\"utzenberger algorithms,
an elementary approach,
{\em Electron. J. Comb.} {\bf 3} (1996), RP R15, 32 p.

\bibitem{L2} M.~A.~A. van Leeuwen,
The Littlewood-Richardson rule, and related combinatorics, in
{\em Interaction of combinatorics and representation theory},
95--145,  {\em MSJ Mem.} {\bf 11}, Math. Soc. Japan, Tokyo, 2001.

\bibitem{Ma} I.~G. Macdonald,
{\em Symmetric functions and Hall polynomials}, second ed.,
Oxford University Press, Oxford, 1995.

\bibitem{P2} I. Pak,
Hook length formula and geometric combinatorics,
{\em S\'emin. Lothar. Comb.} {\bf 46}, B46f (2001), 13 pp.


\bibitem{P1} I. Pak,
Partition Bijections, a Survey,
to appear in {\em Ramanujan Journal},
available from
{\tt http://www-math.mit.edu/} $\widetilde\, $
{\tt pak}.


\bibitem{P3} I. Pak,
Partition Identities and Geometric Bijections, {\em Proc. A.M.S.},
to appear, 2003.


\bibitem{PV} I. Pak and E. Vallejo,
Combinatorics and geometry of Littlewood-Richardson cones,
to appear in {\em Europ. J. Comb.}, available from
{\tt http://www-math.mit.edu/} $\widetilde\,$
{\tt pak}.

\bibitem{Pap} C.~H. Papadimitriou,
{\em Computational Complexity},  Addison-Wesley,
Amsterdam, 1994.

\bibitem{RR} D.~P. Robbins and H. Rumsey,
Determinants and alternating-sign matrices.
{\em Adv. Math.} \, {\bf 62} (1986), 169--184.

\bibitem{Ro} G.~de~B. Robinson,
On the representations of the symmetric group,
{\em Amer. J. Math.} \, {\bf 60} (1938), 745--760.

\bibitem{Sa}
B.~E. Sagan,
{\em The symmetric group. Representations, combinatorial algorithms,
and symmetric functions}, second ed.,
Graduate Texts in Mathematics {\bf 203}, Springer, New York, 2001.

\bibitem{Sch}
C. Schensted,
Longest increasing and decreasing subsequences,
{\em Canad. J. of Math.} {\bf 13} (1961), 179--191.

\bibitem{Sc1}
M.-P. Sch\"utzenberger,
Quelques remarques sur une construction de Schensted,
{\em Math. Scand.} \, {\bf 12} (1963), 117--128.

\bibitem{Sc}
M.-P. Sch\"utzenberger,
Promotion des morphismes d'ensembles ordonn\'es,
{\em Discrete Math.} \, {\bf 2} (1972), 73--94.

\bibitem{Sc3}
M.-P. Sch\"utzenberger,
La correspondance de Robinson,
in {\em Combinatoire et R\'e\-pre\-sentation du Groupe Sym\'etrique},
Lecture Notes in Math. {\bf 579}, Springer, 1977, 59--135.

\bibitem{Sh}
A. Shamir, Factoring Numbers in $O(\log n)$ Arithmetic Steps,
{\em Inf. Proc. Letters} {\bf 8} (1979), 28--31.

\bibitem{Sp}
 D.~E. Speyer,
Perfect Matchings and the Octahedron Recurrence,
arXiv:math.CO/ 0402452.


\bibitem{St} R.~P. Stanley,
{\em Enumerative Combinatorics}, vol. 2,
Cambridge Studies in Advanced Mathematics {\bf 62},
Cambridge University Press, 1999.


\bibitem{stem}
J.~R. Stembridge,
Computational aspects of root systems,
Coxeter groups, and Weyl characters, in
{\em Interaction of combinatorics and representation theory},
1--38 ,  {\em MSJ Mem.} {\bf 11}, Math. Soc. Japan, Tokyo, 2001.


\bibitem{Th}
G.~P. Thomas,
On a construction of Sch\"utzenberger,
{\em Discrete Math.} {\bf 17} (1977), 107--118.

\bibitem{Thu} D. Thurston,  From Dominoes to Hexagons,
arXiv:math.CO/0405482.


\end{thebibliography}
}
\end{document}